\documentclass[12pt]{article}

\usepackage{amssymb}
\usepackage{amsthm}
\usepackage{amsmath}
\usepackage[T1]{fontenc}
\usepackage[utf8]{inputenc}
\usepackage[english]{babel}
\usepackage[pdftex]{graphicx}      
\usepackage{mathrsfs}
\usepackage[sort&compress,numbers]{natbib}

\usepackage{indentfirst}              
\usepackage{makeidx}                    
\usepackage[nottoc]{tocbibind}         
\usepackage{courier}                   
\usepackage{type1cm}                    
\usepackage{listings}                
\usepackage{titletoc}
\usepackage[all,cmtip]{xy}

\usepackage{fancyhdr}
\newtheorem{defi}{\bf Definition}[section]
\newtheorem{theo}[defi]{\bf Theorem}

\newtheorem{coro}[defi]{\bf Corollary}
\newtheorem{pro}[defi]{\bf Proposition}
\newtheorem{lem}[defi]{\bf Lemma}
\newtheorem{ex}[defi]{\bf Example}

\newenvironment{dem}{\noindent
\textnormal{\textbf{Proof:}}}{\begin{flushright}$\mathbf{Q.E.D.}$\end{flushright}}

\theoremstyle{remark}

\newcommand{\algA}{\mathcal{A}}

\def\bee{\begin{eqnarray}}
\def\bes{\begin{eqnarray*}}
\def\eee{\end{eqnarray}}
\def\ees{\end{eqnarray*}}

\title{On a problem by Nathan Jacobson for Malcev algebras}

\author{
{Victor L\'opez Sol\'{\i}s
}\\
{\small Departamento Académico de Ciencias Básicas y Afines}\\
{\small Universidad Nacional de Barranca,}\\
{\small Jr. Toribio de Luzuriaga N°\,376, Barranca, 15169},\\
{\small Lima. Peru}\\
{\small vlopez@unab.edu.pe}
}
\date{\quad}

\begin{document}
\maketitle\vspace{-1.5cm}


\begin{abstract}
%
%
In this paper we solve a problem for a certain class of Malcev algebras, which is an analogous of
an old problem posed by Nathan Jacobson for alternative algebras. Specifically, we prove a coordinatization theorem for a class of Malcev algebras $\mathcal{M}$ containing the 3-dimensional simple Lie algebra $\mathfrak{s l}_{2}(\mathbb{F})$ such that $m\,\mathfrak{s l}_{2}(\mathbb{F})\neq 0$ for any $0\neq m\in\mathcal{M}.$ We drop the last condition and we describe the structure of the same class of Malcev algebras $\mathcal{M}$ that contains $\mathfrak{s l}_{2}(\mathbb{F})$. 
\end{abstract}

{\parindent= 4em \small  \sl Keywords: Malcev algebra, non-Lie Malcev module, Kronecker factorization theorem, Plücker relations.\\
MSC classification (2020): Primary 17D10, secondary 17B22, 17B60, 17B70.}


\section{Introduction}\label{s1}

\subsection{Overview}\label{subsection 1}

A \textit{Malcev algebra} is a $\mathbb{F}$-vector space $\mathcal{M}$ with a bilinear binary operation $(x,y)\mapsto xy $ satisfying the following identities:
\begin{equation}\label{e}
x^2=0, ~~ J(x,y,xz)=J(x,y,z)x, 
\end{equation}
where $J(x,y,z)=(xy)z+(yz)x+(zx)y$ is the Jacobian of the elements $x,y,z\in \mathcal{M}$.

Lie algebras fall into the variety of Malcev algebras because the Jacobian of any three elements vanish. The tangent space $T(L)$ of an analytic Moufang loop $L$ is another example of Malcev algebra. 
Let $\algA$ be an alternative algebra, if we introduce a new product by means of the commutator $[x,y]=xy-yx$ into $\algA$, we obtain a new algebra that will be denoted by $\algA^{(-)}$. It is easy to verify that the algebra $\algA^{(-)}$ is a Malcev algebra. All Malcev algebras obtained in this form are called \textit{special}. 
The classic example of a non-Lie Malcev algebra 
is formed by traceless elements of the Cayley-Dickson algebra with the commutator. This 7-dimensional exceptional Malcev algebra is denoted by $\mathbb{M}$ and is one of most cornerstone examples.

The description of the structure of algebras and superalgebras that contain certain finite-dimensional algebras and superalgebras has a rich history, and important applications in representation theory and category theory (for example, see \cite{J2, LDS1, LSS}). The classical Wedderburn Theorem says that if a unital associative algebra $\algA$ contains a central simple subalgebra of finite dimension $\mathcal{B}$ with the same identity element, then $\algA$ is isomorphic to a Kronecker product $S\otimes_{\mathbb{F}} \mathcal{B}$, where $S$ is the subalgebra of the elements that commute with each $b\in \mathcal{B}$. In particular, if $\algA$ contains $M_{n}(\mathbb{F})$ as subalgebra with the same identity element, we have $\algA\cong M_{n}(S)$ “coordinated” by $S$. Kaplansky in Theorem $2$ of \cite{K} generalized the Wedderburn result to the alternative algebras $\algA$ and the split Cayley algebra $\mathcal{B}$. Jacobson in Theorem $1$ of \cite{J1} gave a new proof of the result of Kaplansky using his classification of completely reducible alternative bimodules over a field of characteristic different of $2$ and finally Victor López-Solís in \cite{LSS} proved that this result is valid for any characteristic. Using this result, Jacobson \cite{J1} proved a Kronecker Factorization Theorem for Jordan algebras that contain the Albert algebra with the same identity element. The statements of this type are usually called \textit{Kronecker factorization theorems}. In the case of right alternative algebras, S. Pchelintsev, O. Shashkov and I. Shestakov \cite{PSS} proved that every unital right alternative bimodule over a Cayley algebra (over an algebraically closed field of characteristic not 2) is alternative and they used that result to prove a coordinatization theorem for unital right alternative algebras containing a Cayley subalgebra with the same unit.

In the case of superalgebras, M. López-Díaz and I. Shestakov \cite{LDS2, LDS1} studied the representations of simple alternative and exceptional Jordan superalgebras in characteristic $3$ and through these representations, they obtained some analogues of the Kronecker Factorization Theorem for these superalgebras. Also, V. López-Solís in \cite{LSS} obtained analogues of the Kronecker Factorization Theorem for some simple alternative superalgebras. Similarly, C. Martinez and E. Zelmanov \cite{CE} obtained a Kronecker Factorization Theorem for the 10-dimensional exceptional Kac superalgebra $K_{10}$. Moreover, A. Pozhidaev and I. Shestakov \cite{PS} studied the representations of simple finite-dimensional noncommutative Jordan superalgebras and proved some analogues of the Kronecker factorization theorem for such superalgebras (see also \cite{PY}). 

In regards to Lie theory, \cite{BM} S. Berman and R. Moody initiated the study of Lie algebras graded by the root system of characteristic 0 and obtained coordinatization theorems for Lie algebras graded by a simply-laced finite root system of rank $n\geq 2$. Moreover, in \cite{BZ} G. Benkart and E. Zelmanov determinated some coordinatization theorems for Lie algebras graded by doubly-laced finite root systems (see also Neher \cite{N}).

Carlson in \cite{C} studied  Malcev modules for 3-dimensional simple Lie algebras and for the 7-dimensional exceptional Malcev algebra $\mathbb{M}$. She proved that the only non-Lie modules that appear for algebraically closed fields of characteristic $\neq 2,3$ are a 2-dimensional module for the 3-dimensional simple Lie algebra $\mathfrak{s l}_{2}(\mathbb{F})$ and the regular module for $\mathbb{M}$. In \cite{CA, EL} it was proved that for fields of characteristic $\neq 2,3$ any finite-dimensional module over a semisimple algebra which is the direct sum of classical simple Lie algebras and Malcev algebras of type $\mathbb{M}$ splits into a direct sum of a Lie module and a direct sum of irreducible non-Lie modules of the two types mentioned above. 
Elduque and I. Shestakov in \cite{ES} determined the  irreducible modules over Malcev algebras without any restriction on the dimension. The main results obtained by them assert that, essentially, the non-Lie modules that appear in the finite-dimensional case are the only one in general.
So using these results, V. López-Solís in \cite{LSS0} proved an analogue of the Kronecker factorization theorem for Malcev (super)algebras that contain (in the even part) $\mathbb{M}$ over fields of characteristic $\neq 2,3$. 
Thus following the study of the structure of (super)algebras that contain (in the even part) finite-dimensional simple algebras, appear the question: what about the structure of Malcev (super)algebras that contain (in the even part) $\mathfrak{s l}_{2}(\mathbb{F})$?.
This question is an analogous of a problem posed by Nathan Jacobson for alternative algebras \cite{J1}.

Analogous of the above problem have been solved for some non-associative algebras. For unital alternative algebras containing a split quaternion algebra with the same identity element has been solved by V. López-Solís and I. Shestakov in \cite{LSS1} (see also \cite{J1}). Similarly, for unital Jordan algebras containing the Jordan algebra of symmetric $2\times 2$ matrices
has been solved by J. Laliena, V. López-Solís and I. Shestakov in \cite{JVI}.


\subsection{Main result}\label{subsection 2}

In this paper, we consider an important subvariety of Malcev algebras. Certainly, Filippov in \cite[(6)]{F2} introduced a very interesting identity and considered (see \cite{F4,F,F1,F3}) the subvariety $\mathcal{H}$ of those Malcev algebras satisfying this identity. He proved in \cite{F} that, in a sense, the subvariety $\mathcal{H}$ plays the same role as the variety of Lie algebras in the structure theory of Malcev algebras. Moreover, the central closures of the prime algebras in $\mathcal{H}$ and the simple algebras in $\mathcal{H}$ are simple algebras of the two types mentioned in the previous subsection (see \cite{ES1}).

The aim of this paper is to solve the proposed problem in subsection \ref{subsection 1} for Malcev algebras in the subvariety $\mathcal{H}$, this is, we determine the structure of the Malcev algebras in $\mathcal{H}$ that contain $\mathfrak{s l}_{2}(\mathbb{F})$. 

The paper is organized as follows: In Section \ref{s2}, we provided some definitions about Malcev algebras, their representations and, graded algebras and modules by the root system. In Section \ref{s3}, we describe the structure of the Malcev algebras $\mathcal{M}$ in $\mathcal{H}$ that contain $\mathfrak{s l}_{2}(\mathbb{F})$.  

Essentially, the coordinatization theorem \ref{principal} involves two ingredients: a Malcev algebra $\mathcal{M}$ in $\mathcal{H}$  containing $\mathfrak{s l}_{2}(\mathbb{F})$ such that $m\mathfrak{s l}_{2}(\mathbb{F})\neq 0$ for any $0\neq m\in\mathcal{M}$ is "coordinated" by a commutative associative algebra $\mathcal{B}$ and by a commutative $\mathcal{B}$-bimodule $V$, on which a skew-symmetric form is defined with values in $\mathcal{B}$, satisfying Plücker relations. More exactly,
\begin{center}
    $\mathcal{M}=\mathfrak{sl}_2(\mathcal{B})\oplus V^2,$
\end{center}
with a properly defined multiplication. Moreover, we drop the  condition "$m\mathfrak{s l}_{2}(\mathbb{F})\neq 0$ for any $0\neq m\in\mathcal{M}$" and we describe the structure of the same class of Malcev algebras $\mathcal{M}$ containing $\mathfrak{s l}_{2}(\mathbb{F})$. 




\section{Preliminaries}\label{s2}

In this section we provide background material that will be used along the paper.

Throughout this paper $\mathbb{F}$ will be a field of characteristic different of 2 and 3. So by $(\ref{e})$, a Malcev algebra $\mathcal{M}$ is an anticommutative algebra that satisfies the Malcev identity
\begin{equation}\label{ee1}
    (xz)(yt)=((xy)z)t+((yz)t)x+((zt)x)y+((tx)y)z.
\end{equation}
Here and in what follows we omit parentheses in left-normed products: $xyzt =
((xy)z)t,$ $xy\cdot zt=(xy)(zt),$ $x\cdot yz=x(yz),$ etc.
The following functions play an important role in
the theory of Malcev algebras:
\begin{equation*}
\begin{split}
 J(x,y,z)&=xyz+yzx+zxy,\mbox{\textit{the Jacobian} of}~x,y,z,\\
\left[x,y,z\right]&=xyz+x\cdot yz,~\mbox{\textit{the antiassociator}},
\end{split}
\end{equation*}
\begin{equation*}
\begin{split}
\{x,y,z\}&= xyz-xz y+2x\cdot yz=J(x,y,z)+3x\cdot yz,\\
h(y,z,t,x,u)&=\{yz,t,u\}x+\{yz,t,x\}u+\{yx,z,u\}t+\{yu,z,x\}t,\\
p(x,y,z,t)&=-\{zt,x,y\}-\{yt,z,x\}+\{xt,y,z\}
\end{split}
\end{equation*}
and the operator 
$$x\alpha(y,z,t)=p(x,y,z,t).$$

The function $p$ is skew-symmetric in $y, z, t.$ Indeed, the skew symmetry of $p$ in $y$ and $z$ follows from the fact that all Malcev algebras are anticommutative. This implies that $\alpha$ is skew-symmetric on its arguments.

Note that the operator $\alpha$ was defined by V. Filippov in \cite{F1} and was essential in \cite{E, ES, F3, LSS0, Sh}.

\subsection{Malcev Modules}\label{subsection 2.1}

Recall that $\mathfrak{s l}_{2}(\mathbb{F})$ is the Lie algebra of traceless $2 \times 2$ matrices over $\mathbb{F}$ with the commutator, and with the following basis $\{E,H,F\}$ such that
\begin{equation}\label{Lie}
    EH=E,\,\,\,FH=-F,\,\,\,EF=\dfrac{1}{2}H.
\end{equation}
This algebra is simple and it is one of the main objects of our paper.

Given a Malcev algebra $\mathcal{M}$ and a vector space $V$ over $\mathbb{F},$ $V$ is said to be a \textit{module} for $\mathcal{M}$ if there exists a $\mathbb{F}$-linear map $\rho:\mathcal{M}\longrightarrow \mbox{End}_{\mathbb{F}}(V)(x\mapsto  \rho_{x})$, such that the split extension $\mathscr{E}=\mathcal{M}\oplus V$, with multiplication given by
\begin{equation}\label{eM}
    (x+v)\cdot(y+w):=xy+(v\rho_{y}-w\rho_{x})
\end{equation}
for $x,y\in \mathcal{M}$, $v,w\in V,$
is a Malcev algebra. In this case, $\rho$ is called a \textit{representation} of $\mathcal{M}$.

The module $V$ is said to be \textit{irreducible} if $\rho\neq 0$ and it does not contain any proper submodule. $V$ is said to be Lie module in case $J(V,\mathcal{M},\mathcal{M})=0$ in $\mathscr{E}$. Also, $V$ is said to be \textit{almost faithful} if $\mbox{ker}\rho$ does not contain any ideal of $\mathcal{M}$ and faithful if $\mbox{ker}\rho=0.$ 

A \textit{$\mathcal{M}$-adjoint module}, $\overline{\mathcal{M}}$, for algebra $\mathcal{M}$, is defined on the vector space $\mathcal{M}$ with the action of $\mathcal{M}$ coinciding with the multiplication in $\mathcal{M}$. So we have the $\mathfrak{sl}_{2}(\mathbb{F})$-adjoint module, $\overline{\mathfrak{sl}_{2}(\mathbb{F})}.$

\begin{ex}\label{ex1}
The Lie algebra $\mathfrak{s l}_{2}(\mathbb{F})$ has a 2-dimensional non-Lie Malcev module with basis $\{u,v\}$ such that
\begin{center}
    $uH=u,\,\,vH=-v,\,\, uE=v,\,\, uF=0,\,\, vE=0,\,\, vF=-u$.
\end{center}
In characteristic 3, this is a Lie module. This $\mathfrak{s l}_{2}(\mathbb{F})$-module will be called a module of type $V_{2}$.
\end{ex}

Let $V$ be a module for the Malcev algebra $\mathcal{M}$ and let $\mathscr{E}=\mathcal{M}\oplus V$ be the corresponding split extension. Let us consider
$$\Gamma=\Gamma(\mathscr{E})=\{\alpha\in\mbox{End}_{\mathbb{F}}(\mathscr{E}):(xy)\alpha=x(y\alpha)=(x\alpha)y\,\,\,\forall\,x,y\in \mathscr{E}\},\,\,\mbox{\textit{the centroid} of}\,\,\mathscr{E},$$
$$Z=Z(\mathscr{E})=\{\alpha\in\Gamma:~V\alpha\subseteq V~\mbox{and}~\mathcal{M}\alpha\subseteq \mathcal{M}\},$$
$$K=K(V)=\{\varphi\in\mbox{End}_{\mathbb{F}}(V):\, \varphi\rho_x=\rho_x \varphi \,\,\forall\,x\in \mathcal{M}\},\,\,\mbox{\textit{the centralizer} of}\,\,V.$$

The following proposition was proved in (\cite{ES}, Proposition $4$) to superalgebras case. Essentially provides some basic properties of the subalgebras $Z$ and $K$ when we consider irreducible almost faithful modules. 

\begin{pro}\label{p1}
Assume that $V$ is an irreducible almost faithful module for $\mathcal{M}$. Then,
\begin{description}
\item[(i)] $Z$ is an integral domain which acts without torsion on $\mathcal{M}$.
\item[(ii)] $K$ is a skew field and any nonzero element in $K$ acts bijectively on $V$.
\item[(iii)] The restriction homomorphism $\phi:Z\longrightarrow K(\alpha\longmapsto\alpha|_V)$ is one-to-one.
\end{description}
\end{pro}

\subsection{Graded algebras and modules by the root system} 

Let $G$ be any group with group unit 1. 
If $A$ is an algebra, we say $A$ is a $G$-\textit{graded algebra} if there are subspaces $A_g$, for each $g\in G$, such that
\begin{center}
    $A=\bigoplus_{g\in G} A_g$ \textup{and for each} $g,h\in G,$ $A_g A_h\subseteq A_{gh}.$
\end{center}
If $0\neq a \in A_g,$ we say $a$ is \textit{homogeneous of $G$-degree} $g$ and we write $\textup{deg}(a) = g.$

We know that the Lie algebra $\mathfrak{sl}_{2}(\mathbb{F})$ has the triangular decomposition $\mathfrak{sl}_{2}(\mathbb{F})=L_{-1}\oplus L_{0}\oplus L_{1}$ determined by the Cartan subalgebra $L_0=\langle H\rangle,$ where $L_1=\langle E\rangle$ and  $L_{-1}=\langle F\rangle.$ This decomposition is the so-called \textit{the Cartan grading} of $\mathfrak{sl}_{2}(\mathbb{F}),$ from where we conclude that
\begin{center}
    $\mathfrak{sl}_{2}(\mathbb{F})=\bigoplus_{i\in\mathbf{Z}_{3}}L_{i}$
\end{center}
since $-1=2$ in $\mathbf{Z}_{3}.$ Thus $\mathfrak{sl}_{2}(\mathbb{F})$ is a $\mathbf{Z}_{3}$-graded Lie algebra.

\begin{defi}\label{defig}
A Lie algebra $L$ over $F$ is \textit{graded by the root system} $\Phi$, or
$\Phi$-graded, if:
\begin{center}
    \begin{enumerate}
        \item $L$ contains as a subalgebra a finite-dimensional simple Lie algebra $\mathfrak{g}= \mathfrak{h} \oplus (\bigoplus_{\alpha\in \Phi} \mathfrak{g}_{\alpha})$ whose root system is $\Phi$, relative to a Cartan subalgebra $\mathfrak{h} = \mathfrak{g}_{0}$;
        \item  $L = \bigoplus_{\alpha\in \Phi\cup \{0\}} L_{\alpha}$, where $L_\alpha= \{X \in L \,\,| \,\,[X, H] = \alpha(H)X\,\, \textup{for all}\,\, H \in h\}$; and 
        \item $L_0 = \sum_{\alpha\in\Phi}[L_\alpha, L_{-\alpha}]$.
    \end{enumerate}
\end{center}
\end{defi}
The subalgebra $\mathfrak{g}$ is said to be a \textit{grading subalgebra} of L.

It follows from the Definition \ref{defig} that $\mathfrak{sl}_{2}(\mathbb{F})$ is graded by the root system $\Phi$, where $\Phi=\{-1,1\}.$

\begin{defi}\label{defig1} 
A  Lie module $V$ over a G-graded Lie algebra $L$ is called
$G$-graded if $V$ is a $G$-graded vector space and
\begin{center}
    $L_{g}V_{h} \subseteq V_{gh}$
\end{center}
for all $g, h \in G$.
\end{defi}

It is well known that the $\mathfrak{sl}_{2}(\mathbb{F})$-adjoint module is an irreducible Lie module. As above, we denote this module by $\overline{\mathfrak{sl}_{2}(\mathbb{F})}$ and, for $a\in \mathfrak{sl}_{2}(\mathbb{F})$, the corresponding element
in $\overline{\mathfrak{sl}_{2}(\mathbb{F})}$ will be denoted by $\overline{a}$. Therefore, the action of $\mathfrak{sl}_{2}(\mathbb{F})$ in $\overline{\mathfrak{sl}_{2}(\mathbb{F})}$ satisfies 
\begin{align*}
E\overline{E}&=0,           &  E\overline{H}&=\overline{E},              &  E\overline{F}&=\frac{1}{2}\overline{H},\\
H\overline{E}&=-\overline{E},   &  H\overline{H}&=0,         & H\overline{F}&=\overline{F},\\
F\overline{E}&=-\frac{1}{2}\overline{H},       &  F\overline{H}&=-\overline{F},   &  F\overline{F}&=0.
\end{align*}

It is clear that $\overline{\mathfrak{sl}_{2}(\mathbb{F})}$ is a $\mathbf{Z}_{3}$-graded module; $\overline{\mathfrak{sl}_{2}(\mathbb{F})}=V_{-1}\oplus V_{0}\oplus V_{1}$, where $V_{1}=\langle \overline{E}\rangle,$ $V_{0}=\langle \overline{H}\rangle$ and $V_{-1}=\langle \overline{F}\rangle$.

\section{Malcev algebras that contain $\mathfrak{sl}_{2}(\mathbb{F})$}\label{s3}


Let $H(\mathcal{M})$ be the subspace spanned by the elements $h(y,z,t,u,x)$; $H(\mathcal{M})$ is an ideal of the Malcev algebra $\mathcal{M}$ (see \cite{F}). 
The subvariety $\mathcal{H}$ (over $\mathbb{F}$) is defined as the class of Malcev algebras that satisfy the identity $h(y,z,t,u,x)=0$, that is, $H(\mathcal{M})=0$. It is proved in \cite{F1} that a Malcev algebra $\mathcal{M}$ belongs of $\mathcal{H}$ if and only if $\alpha(\mathcal{M},\mathcal{M},\mathcal{M})$ (the linear span of the operators $\alpha(x,y,z)$) is contained in the $\Gamma(\mathcal{M})$ and in this case $\alpha(\mathcal{M},\mathcal{M},\mathcal{M})$ is an ideal of the centroid.

The most important finite-dimensional simple algebras in the variety $\mathcal{H}$ are the 3-dimensional Lie algebra $\mathfrak{sl}_{2}(\mathbb{F})$ and the 7-dimensional exceptional Malcev algebra $\mathbb{M}$ over its centroid $\Gamma$, which is a field (see \cite{E}, Theorems 2 and 4).




Let $\mathcal{M}$ be a Malcev algebra in $\mathcal{H}$ and let $V$ be a module for $\mathcal{M}$, $V$ is said to be a module for $\mathcal{M}$ in the variety $\mathcal{H}$ or $\mathcal{H}$-\textit{module} if there exists a $\mathbb{F}$-linear map $\rho:\mathcal{M}\longrightarrow \mbox{End}_{\mathbb{F}}(V)(x\mapsto  \rho_{x})$ such that the split extension $\mathscr{E}=\mathcal{M}\oplus V$ is again a Malcev algebra in $\mathcal{H}$ with a multiplication $(\ref{eM})$.

Analogously, we can define the subvariety $\mathcal{H}$ to superalgebras case, which has the same notation as in the algebras. The following results were proved by A. Elduque (\cite{E}, Theorem 5 and Corollary 6) to Malcev superalgebras in $\mathcal{H}$ and its proofs can be easily adapted to algebras case.

\begin{theo}\label{t1}
Let $\mathcal{M}$ be a simple algebra in the variety $\mathcal{H}$ and let $V$ be a $\mathcal{H}$-module for $\mathcal{M}$. Then $V$ is completely reducible.
\end{theo}

As a consequence of the Theorem $\ref{t1}$, we can deduce:
\begin{coro}\label{corol2}
Let $\mathcal{M}$ be a simple algebra in $\mathcal{H}$, $V$ an $\mathcal{H}$-module for $\mathcal{M}$, and $\Gamma$ the centroid of $\mathcal{M}$. Then $V=N_{V} \oplus J_{V}$ where the submodules $N_{V}$ and $J_{V}$ are given by $N_{V}=\{v \in V: J(v, \mathcal{M}, \mathcal{M})=0\}$ and $J_{V}=J(V, \mathcal{M}, \mathcal{M})$. Moreover,

\begin{description}
\item[(i)] If $\operatorname{dim}_{\Gamma} \mathcal{M}=3$, then $N_{V}=\textup{Ann}_{V}\mathcal{M} \oplus \widehat{N}_{V}$, with $\textup{Ann}_{V}\mathcal{M}=\{v \in V: v \mathcal{M}=0\}$, $\widehat{N}_{V}=V \mathcal{M} \cap N_{V} .$ Besides, $\widehat{N}_{V}$ is a direct sum of copies of the adjoint module for $\mathcal{M}$ and $J_{V}$ is a direct sum of copies of the unique irreducible non-Lie module for $\mathcal{M}$.
\item[(ii)] If $\operatorname{dim}_{\Gamma} \mathcal{M}=7$, then $N_{V} \mathcal{M}=0$ and $J_{V}$ is a direct sum of copies of the adjoint module for $\mathcal{M}$.
\end{description}
\end{coro}
Moreover, if $\mathcal{M}=\mathfrak{sl}_{2}(\mathbb{F})$ in the Corollary \ref{corol2}$(i)$, then I. Shestakov (\cite{Sh}, Lemma 8) proved that
$J_{V}=\{v\in V:\{v,\mathcal{M},\mathcal{M}\}=0\}$, that is, $J_{V}$
satisfies the condition
\begin{center}
    $vab-vba+2v\cdot ab=0.$
\end{center}
for any $a,b\in \mathcal{M}$ and $v\in J_{V}$.

Now if $V$ is a Malcev module for $\mathfrak{sl}_{2}(\mathbb{F})$ not necessarily in the variety $\mathcal{H}$, it should be remarked that $N_{V}$ need not be completely reducible (see, for instance \cite{GO}). This is the reason why we  consider Malcev algebras in $\mathcal{H}$ because by Corollary \ref{corol2} $(i)$ $\widehat{N}_{V}$ is completely reducible and is a direct sum of copies of the adjoint module for $\mathfrak{sl}_{2}(\mathbb{F})$.

But in the finite-dimensional case the situation is different because it was proved that any finite-dimensional Malcev $\mathfrak{sl}_{2}(\mathbb{F})$-module splits into a direct sum of a Lie module completely reducible and modules of type $V_{2}$. The class of irreducible Lie modules is comprised of modules of type $L_{m}$, where $m$ is a non-negative integer and $L_{m}$ has the basis $\{v_{0},v_{1},\dots, v_{m}\}$ with multiplication table
$$
\begin{gathered}
v_{i} H=\frac{1}{2}(m-2 i)v_i, \quad i=0, \ldots, m \\
v_{i} E=\frac{1}{2}v_{i+1}, \quad i=0, \ldots, m-1, \quad v_{m} E=0 \\
v_{i} F=\frac{1}{2}(m i-(i-1) i) v_{i-1}, \quad i=1, \ldots, m, \quad v_{0} F=0.
\end{gathered}
$$

In \cite{LSS0} was described the structure of Malcev algebras containing the 7-dimensional exceptional Malcev algebra $\mathbb{M}$. In this paper we describe the structure of Malcev algebras $\mathcal{M}$ in the variety $\mathcal{H}$ containing $\mathfrak{sl}_{2}(\mathbb{F})$.

\subsection{The case $mL\neq 0$ for any $0\neq m\in\mathcal{M}$}\label{3.1}

Let $\mathcal{M}$ be a Malcev algebra in the variety $\mathcal{H}$ that contains a subalgebra $L=\mathfrak{sl}_{2}(\mathbb{F})$ such that $mL\neq 0$ for any $0\neq m\in\mathcal{M}.$ 
Then we can consider $\mathcal{M}$ as a Malcev $\mathcal{H}$-module for $L$ and by Corollary $\ref{corol2}$$(i)$, $\textup{Ann}_{\mathcal{M}}L=0$, $\widehat{N}_{\mathcal{M}}=N_{\mathcal{M}}$ and
\begin{center}
$\mathcal{M}=N_{\mathcal{M}}\oplus J_{\mathcal{M}}$,
\end{center}
where $N_{\mathcal{M}}=\displaystyle\sum_{i}\oplus\overline{\mathfrak{sl}_{2}(\mathbb{F})}_{i}$ with $\overline{\mathfrak{sl}_{2}(\mathbb{F})}_{i}\cong \overline{\mathfrak{sl}_{2}(\mathbb{F})}$ as $\mathfrak{sl}_{2}(\mathbb{F})$-modules and 
$J_{\mathcal{M}}=\sum_{i}\oplus V_{2i},$ with $V_{2i}\cong V_2$ as $\mathfrak{sl}_{2}(\mathbb{F})$-modules, where $V_{2i}$ is a 2-dimensional non-Lie Malcev module for $\mathfrak{sl}_{2}(\mathbb{F})$ which has a basis $\{u_{i},v_{i}\}$ (see Example \ref{ex1}) satisfying
\begin{equation}\label{e5}
u_{i}H=u_{i},\,\,v_{i}H=-v_{i},\,\, u_{i}E=v_{i},\,\, u_{i}F=0,\,\, v_{i}E=0,\,\, v_{i}F=-u_{i}.
\end{equation} 

The expression $N_{\mathcal{M}}=\displaystyle\sum_{i}\oplus\overline{\mathfrak{sl}_{2}(\mathbb{F})}_{i}$ allows us to conclude that $N_{\mathcal{M}}$ is a $\mathbf{Z}_3$-graded module, that is, 
\begin{equation*}
    N_{\mathcal{M}}=V_{-1}\oplus V_{0}\oplus V_{1},
\end{equation*}
where $V_{1}=\langle \overline{E}_{i}\rangle,$ $V_{0}=\langle \overline{H}_{i}\rangle$ and $V_{-1}=\langle \overline{F}_{i}\rangle.$ 

Moreover by the proof of the Theorem 5 of \cite{E}, we have that
\begin{center}    $N_{\mathcal{M}}=L\alpha(\mathcal{M},L,L)=\displaystyle\sum_{0\neq\alpha\in\alpha(\mathcal{M},L,L)}\oplus L\alpha$
\end{center}
is the sum of copies of $\overline{L}$ and $\alpha(\mathcal{M},L,L)$ is a subspace (the linear span of the operators $\alpha(x,y,z),$ where $x\in\mathcal{M}$ and $y,z\in L$) of the ideal $\alpha(\mathcal{M},\mathcal{M},\mathcal{M})$ of $\Gamma(\mathcal{M})$. 

In Proposition \ref{c1} we prove a  coordinatization theorem for $N_{\mathcal{M}}$, this means that $N_{\mathcal{M}}$ is isomorphic to $\mathfrak{sl}_{2}(U)$ as algebras, where $U$ is a certain associative commutative algebra.

\subsubsection{$\mathbf{Z}_2$-grading $\mathcal{M}=N_{\mathcal{M}}\oplus J_{\mathcal{M}}$}


\begin{lem}\label{l2}
$N_{\mathcal{M}}$ is a subalgebra of $\mathcal{M}$ and $N_{\mathcal{M}}J_{\mathcal{M}}\subseteq J_{\mathcal{M}}$.
\end{lem}
\begin{dem}
We have $N_{\mathcal{M}}=\{m\in  \mathcal{M}:J(m,L,L)=0\}$ and $J_{\mathcal{M}}=\{m\in \mathcal{M}:\{m,L,L\}=0\}$. Then for any $\alpha, \beta\in \alpha(\mathcal{M},L,L)$
\begin{center}
    $J((L\alpha)(L\beta),L,L)\subseteq J(L^{2}\alpha\beta,L,L) \subseteq J(L,L,L)\alpha\beta=0$,
\end{center}
so $(L\alpha)(L\beta)\subseteq N_{\mathcal{M}}$ and $N_{\mathcal{M}}N_{\mathcal{M}}\subseteq N_{\mathcal{M}}$. Consequently $N_{\mathcal{M}}$ is a subalgebra of $\mathcal{M}$.

Furthermore 
\begin{center}
    $\{(L\alpha)J_{\mathcal{M}},L,L\}\subseteq \{(L J_{\mathcal{M}})\alpha,L,L\}\subseteq \{J_{\mathcal{M}},L,L\}\alpha=0$,
\end{center}
then $(L\alpha)J_{\mathcal{M}}\subseteq J_{\mathcal{M}}$ and $N_{\mathcal{M}}J_{\mathcal{M}}\subseteq J_{\mathcal{M}}$.

Thus the lemma is proved.
\end{dem}

\begin{pro}\label{c1}
The subalgebra $N_{\mathcal{M}}$ is isomorphic to
$L\otimes_{\mathbb{F}} U$, where $U$ is a unital associative commutative algebra.
\end{pro}
\begin{dem} 
We have $N_{\mathcal{M}}=L\alpha(\mathcal{M},L,L)$.
Denote $U=\alpha(\mathcal{M},L,L),$ and we will show that $U$ is the desired subalgebra of $\Gamma(\mathcal{M})$. 

We fix arbitrary elements $\alpha, \beta\in U$ and $a,b\in L$. So $a\alpha\beta\in N_{\mathcal{M}}$ because
\begin{center}
    $J(a\alpha\beta,L,L)\subseteq J(a\alpha,L,L)\beta=0$,
\end{center}
then $a\alpha\beta=\sum_i a_i\alpha_i$ for some $a_i\in L$ and $\alpha_i\in U$. Assume that $w=a\alpha\beta-\sum_i a_i\alpha_i=0$,
\begin{equation}\label{e7}
    w=(\delta E+ \lambda F+ \gamma H)\alpha\beta-\sum_{i}(\delta_{i} E+ \lambda_{i} F+ \gamma_{i} H)\alpha_i=0,
\end{equation}
where $a=\delta E+ \lambda F+ \gamma H$, $a_i=\delta_{i} E+ \lambda_{i} F+ \gamma_{i} H$ and $\delta, \lambda, \gamma, \delta_{i}, \lambda_{i}, \gamma_{i}\neq 0$ are elements of $\mathbb{F}$.
Since $\alpha\beta\in \Gamma(\mathcal{M})$ and by the multiplication table of $L$, from $(\ref{e7})$ we get
\begin{center}
    $0=Ew=(\frac{1}{2}\lambda H+ \gamma E)\alpha\beta-\sum_{i}(\frac{1}{2}\lambda_{i} H+ \gamma_{i} E)\alpha_i,$
\end{center}
\begin{equation*}
    0=E\cdot Ew=(\frac{1}{2}\lambda E)\alpha\beta-\sum_{i}(\frac{1}{2}\lambda_{i} E)\alpha_i
\end{equation*}
which implies $E\alpha\beta=E\widetilde{\alpha}$, where $\widetilde{\alpha}=\sum_{i}\lambda_{i}\lambda^{-1}\alpha_i\in U$. Then $E(\alpha\beta-\widetilde{\alpha})=0$. From this it is easy to see that $F(\alpha\beta-\widetilde{\alpha})=0$ and $H(\alpha\beta-\widetilde{\alpha})=0$; thus
\begin{center}
    $L(\alpha\beta-\widetilde{\alpha})=0$
\end{center}
and $\alpha\beta-\widetilde{\alpha}=0$. Therefore $\alpha\beta=\widetilde{\alpha}\in U$; $UU\subseteq U$.

We consider $m\in \mathcal{M}$ and $a\in L$, then
\begin{equation*}
\begin{split}
((m\alpha)\beta)a&=(m\alpha)(a\beta)=((m\alpha)a)\beta=(m(a\alpha))\beta\\
                    &=(m\beta)(a\alpha)=((m\beta)a)\alpha=((m\beta)\alpha)a.
\end{split}
\end{equation*}
If $[\alpha,\beta]=\alpha\beta-\beta\alpha$, we have 
\begin{center}
    $(\mathcal{M}[\alpha,\beta])L=0;$
\end{center}
so $\mathcal{M}[\alpha,\beta]=0$ (because $\textup{Ann}_{\mathcal{M}}L=0$) and $[\alpha,\beta]|_{\mathcal{M}}=0$. In particular, $[\alpha,\beta]|_{W_{i}}=0$ for any irreducible component $W_i$ of $\mathcal{M}$, then by a version for algebras of the Proposition $\ref{p1}$(iii) we get $[\alpha,\beta]=0$ because $\phi:Z\longrightarrow K$($\alpha\longmapsto\alpha|_{W_{i}}$) is one-to-one. Therefore, $[U,U]=0$; then $U$ is a unital associative commutative algebra.

Also,
\begin{equation*}
\begin{split}
(a\alpha)(b\beta)&=(a(b\beta))\alpha=((ab)\beta)\alpha\\
                    &=(ab)(\beta\alpha)=(ab)(\alpha\beta).
\end{split}
\end{equation*}

Let $v=E\alpha_{1}+F\alpha_{2}+H\alpha_{3}=0$, where $\alpha_i\in U$. Then,
$0=Ev=\frac{1}{2}H\alpha_{2}+E\alpha_{3}$ and
$0=E\cdot E v=\frac{1}{2}E\alpha_2$
which implies $E\alpha_2=0$. From this, $F\alpha_2=0$ and $H\alpha_2=0$; so
\begin{center}
    $L\alpha_{2}=0$
\end{center}
and $\alpha_2=0$. Similarly, $\alpha_i=0$ for $i=1,3$. Therefore, $N_{\mathcal{M}}\cong L\otimes_{\mathbb{F}} U$ as algebras.

The proposition is proved.
\end{dem}

\begin{coro}\label{coro1}
    The subalgebra $N_{\mathcal{M}}$ of $\mathcal{M}$ is a Lie algebra. 
\end{coro}

\begin{coro}\label{coro2}
    The Lie algebra $N_{\mathcal{M}}$ of $\mathcal{M}$ is $\mathbf{Z}_{3}$-graded perfect Lie algebra.
\end{coro}


Note that the following result does not require the condition: $mL\neq 0$ for any $0\neq m\in\mathcal{M}.$ Therefore, it has a general figure.

\begin{lem}(\cite{G}, Lemma 9)\label{l3}
Let $V_{2i}$ and $V_{2j}$ be the 2-dimensional non-Lie Malcev modules over $L$ with the standard bases $\{u_{i},v_{i}\}$ and $\{u_{j},v_{j}\}$, respectively, then 
\begin{align*}
u_{i}u_{j}E&=-u_{i}v_{j},           &  u_{i}v_{j}E&=-\frac{1}{2}v_{i}v_{j},              &  v_{i}v_{j}E&=0,\\
u_{i}u_{j}F&=0,       &  u_{i}v_{j}F&=\frac{1}{2}u_{i}u_{j},   &  v_{i}v_{j}F&=u_{i}v_{j},\\
u_{i}u_{j}H&=-u_{i}u_{j},   &  u_{i}v_{j}H&=0,         & v_{i}v_{j}H&=v_{i}v_{j}
\end{align*}
and $u_{i}v_{j}=v_{i}u_{j}$ for all $i,j$. So, $J^{2}_{\mathcal{M}}$ is a Lie $L$-module. In particular, if $i=j$, then $V^{2}_{2i}=0$.
\end{lem}
\begin{dem}
From $(\ref{e5})$ for all $i,j,$ we have 
\begin{equation*}
\begin{split}
-u_{i}v_{j}F&=u_{i}v_{j}\cdot FH\overset{(\ref{ee1})}{=}u_{i}Fv_{j}H+Fv_{j}Hu_{i}+v_{j}Hu_{i}F+Hu_{i}Fv_{j}\\
                    &=u_{j}Hu_{i}-v_{j}u_{i}F-u_{i}Fv_{j}=u_{j}u_{i}-v_{j}u_{i}F;
\end{split}
\end{equation*}
so
\begin{equation}\label{e6}
    u_{i}v_{j}F=\frac{1}{2}u_{i}u_{j}.
\end{equation}
Also,
\begin{equation*}
\begin{split}
-v_{i}v_{j}F&=v_{i}v_{j}\cdot FH\overset{(\ref{ee1})}{=}v_{i}Fv_{j}H+Fv_{j}Hv_{i}+v_{j}Hv_{i}F+Hv_{i}Fv_{j}\\
                    &=-u_{i}v_{j}H+u_{j}Hv_{i}-v_{j}v_{i}F+v_{i}Fv_{j}=-u_{i}v_{j}H+u_{j}v_{i}+v_{i}v_{j}F-u_{i}v_{j}
\end{split}
\end{equation*}
and 
\begin{equation}\label{e8}
    v_{i}v_{j}F=\frac{1}{2}(u_{i}v_{j}H-u_{j}v_{i}+u_{i}v_{j}). 
\end{equation}
Similarly, 
\begin{equation*}
\begin{split}
\frac{1}{2}u_{i}v_{j}H&=u_{i}v_{j}\cdot EF\overset{(\ref{ee1})}{=}u_{i}Ev_{j}F+Ev_{j}Fu_{i}+v_{j}Fu_{i}E+Fu_{i}Ev_{j}\\
                    &=v_{i}v_{j}F-u_{j}u_{i}E\overset{(\ref{e8})}{=}\frac{1}{2}(u_{i}v_{j}H-u_{j}v_{i}+u_{i}v_{j})-u_{j}u_{i}E;
\end{split}
\end{equation*}
hence 
\begin{equation}\label{e9}
    u_{i}u_{j}E=-\frac{1}{2}(u_{i}v_{j}+v_{i}u_{j}).
\end{equation}
Moreover,
\begin{equation*}
\begin{split}
u_{i}u_{j}E&=u_{i}u_{j}\cdot EH\overset{(\ref{ee1})}{=}u_{i}Eu_{j}H+Eu_{j}Hu_{i}+u_{j}Hu_{i}E+Hu_{i}Eu_{j}\\
&=v_{i}u_{j}H-v_{j}Hu_{i}+u_{j}u_{i}E-u_{i}Eu_{j}=v_{i}u_{j}H+v_{j}u_{i}+u_{j}u_{i}E-v_{i}v_{j},
\end{split}
\end{equation*}
then
\begin{center}
    $2u_{i}u_{j}E=v_{i}u_{j}H+v_{j}u_{i}-v_{i}v_{j},$
\end{center}
and using $(\ref{e9})$
\begin{center}
    $-u_{i}v_{j}-v_{i}u_{j}=v_{i}u_{j}H+v_{j}u_{i}-v_{i}v_{j},$
\end{center}
we get
\begin{equation}\label{e10}
    v_{i}u_{j}H=0.
\end{equation}
Then, by $(\ref{e8})$, $v_{i}v_{j}F=\frac{1}{2}(u_{i}v_{j}H-u_{j}v_{i}+u_{i}v_{j})=\frac{1}{2}(-u_{j}v_{i}+u_{i}v_{j}).$ 
Further,
\begin{equation*}
\begin{split}
v_{i}v_{j}F&=u_{i}Ev_{j}F
\overset{(\ref{ee1})}{=}u_{i}v_{j}\cdot EF-Ev_{j}Fu_{i}-v_{j}Fu_{i}E-Fu_{i}Ev_{j}\\
&=\frac{1}{2}u_{i}v_{j}H+u_{j}u_{i}E\overset{(\ref{e10})}{=}-u_{i}u_{j}E,
\end{split}
\end{equation*}
so by $(\ref{e9})$ we have
\begin{equation}\label{e11}
    v_{i}v_{j}F=-u_{i}u_{j}E.
\end{equation}
Besides,
\begin{equation*}
\begin{split}
u_{i}v_{j}E&=u_{i}Hv_{j}E\overset{(\ref{ee1})}{=}u_{i}v_{j}\cdot HE-Hv_{j}Eu_{i}-v_{j}Eu_{i}H-Eu_{i}Hv_{j}\\
                    &=-u_{i}v_{j}E-v_{j}Eu_{i}+v_{i}Hv_{j}\\
                    &=-u_{i}v_{j}E-v_{i}v_{j}
\end{split}
\end{equation*}
and 
\begin{equation}\label{e16}
    u_{i}v_{j}E=-\frac{1}{2}v_{i}v_{j}.
\end{equation}
In the same way
\begin{equation*}
\begin{split}
-u_{i}v_{j}&=v_{i}F\cdot u_{j}E\overset{(\ref{ee1})}{=}v_{i}u_{j}FE+u_{j}FEv_{i}+FEv_{i}u_{j}+Ev_{i}u_{j}F\\
&\overset{(\ref{e6})}{=}-\frac{1}{2}u_{j}u_{i}E-\frac{1}{2}Hv_{i}u_{j}=\frac{1}{2}u_{i}u_{j}E-\frac{1}{2}v_{i}u_{j},
\end{split}
\end{equation*}
then $u_{i}u_{j}E=-2u_{i}v_{j}+v_{i}u_{j}$. So from $(\ref{e9})$, $-\frac{1}{2}(u_{i}v_{j}+v_{i}u_{j})=-2u_{i}v_{j}+v_{i}u_{j}$ which implies 
\begin{equation}\label{e18}
u_{i}v_{j}=v_{i}u_{j}.
\end{equation}
Hence, from $(\ref{e9})$ and $(\ref{e11})$, we get
\begin{equation}\label{ee2}
    v_{i}v_{j}F=-u_{i}u_{j}E=u_{i}v_{j}.
\end{equation}
Furthermore,
\begin{equation*}
\begin{split}
-u_{i}u_{j}F&=v_{i}Fu_{j}F\overset{(\ref{ee1})}{=}v_{i}u_{j}\cdot FF-Fu_{j}Fv_{i}-u_{j}Fv_{i}F-Fv_{i}Fu_{j}\\
&=-u_{i}Fu_{j}=0,
\end{split}
\end{equation*}
then 
\begin{equation}\label{ee20}
    u_{i}u_{j}F=0.
\end{equation}
Similarly, we have 
\begin{equation*}
\begin{split}
-u_{i}u_{j}H&=v_{i}Fu_{j}H\overset{(\ref{ee1})}{=}v_{i}u_{j}\cdot FH-Fu_{j}Hv_{i}-u_{j}Hv_{i}F-Hv_{i}Fu_{j}\\
&=-v_{i}u_{j}F-u_{j}v_{i}F-v_{i}Fu_{j}=u_{i}u_{j}
\end{split}
\end{equation*}
and so
\begin{equation}\label{ee21}
    u_{i}u_{j}H=-u_{i}u_{j}.
\end{equation}
Moreover,
\begin{equation*}
\begin{split}
v_{i}v_{j}H&=u_{i}Ev_{j}H\overset{(\ref{ee1})}{=}u_{i}v_{j}\cdot EH-Ev_{j}Hu_{i}-v_{j}Hu_{i}E-Hu_{i}Ev_{j}\\
&=u_{i}v_{j}E+v_{j}u_{i}E+u_{i}Ev_{j}=v_{i}v_{j},
\end{split}
\end{equation*}
implies
\begin{equation}\label{ee22}
    v_{i}v_{j}H=v_{i}v_{j}.
\end{equation}
Further,
\begin{equation*}
\begin{split}
v_{i}v_{j}E&=u_{i}Ev_{j}E\overset{(\ref{ee1})}{=}u_{i}v_{j}\cdot EE-Ev_{j}Eu_{i}-v_{j}Eu_{i}E-Eu_{i}Ev_{j}\\
&=v_{i}Ev_{j}=0,
\end{split}
\end{equation*}
then
\begin{equation}\label{ee23}
    v_{i}v_{j}E=0.
\end{equation}

If $i=j$, then by $(\ref{ee2})$, $u_{i}v_{i}=-u_{i}u_{i}E=0.$ Thus $V^{2}_{2i}=0.$

We consider the identity \cite{E}
\begin{equation}\label{ee34}
    J(tx,y,z)=tJ(x,y,z)+J(t,y,z)x-2J(t,x,yz)
\end{equation}
which is valid in every Malcev algebra. Then, for any $u,v\in J_{\mathcal{M}}$ and $a,b\in L$, we have $\{u,a,b\}=\{v,a,b\}=0$ and from $(\ref{ee34})$
\begin{equation*}
\begin{split}
J(uv,a,b)&=uJ(v,a,b)+J(u,a,b)v-2J(u,v,ab)\\
&=u(\{v,a,b\}-3v\cdot ab)+(\{u,a,b\}-3u\cdot ab)v\\
&\,\,\,\,\,\,-2(uv\cdot ab+v(ab)\cdot u+(ab)u\cdot v)\\
&= v(ab)\cdot u+(ab)u\cdot v-2uv\cdot ab\\
&= \{ab,u,v\}.
\end{split}
\end{equation*}
Thus, 
\begin{equation}\label{ee24}
    J(uv,a,b)=\{ab,u,v\}.
\end{equation}
Then, for any $i,j$ we obtain
\begin{equation}\label{lie1}
\begin{split}
\{E,u_{i},u_{j}\}&=Eu_{i}u_{j}-Eu_{j}u_{i}+2E\cdot u_{i}u_{j}\\
&\overset{(\ref{ee2})}{=}-v_{i}u_{j}+v_{j}u_{i}+2u_{i}v_{j}\overset{(\ref{e18})}{=}-u_{i}v_{j}-u_{i}v_{j}+2u_{i}v_{j}=0,\\
\{F,u_{i},u_{j}\}&=Fu_{i}u_{j}-Fu_{j}u_{i}+2F\cdot u_{i}u_{j}\overset{(\ref{ee20})}{=}0,\\
\{H,u_{i},u_{j}\}&=Hu_{i}u_{j}-Hu_{j}u_{i}+2H\cdot u_{i}u_{j} \overset{(\ref{ee21})}{=} -u_{i}u_{j}+u_{j}u_{i}+2u_{i}u_{j}=0
\end{split}
\end{equation}
and
\begin{equation}\label{lie2}
\begin{split}
\{E,v_{i},v_{j}\}&=Ev_{i}v_{j}-Ev_{j}v_{i}+2E\cdot v_{i}v_{j}\overset{(\ref{ee23})}{=}0,\\
\{F,v_{i},v_{j}\}&=Fv_{i}v_{j}-Fv_{j}v_{i}+2F\cdot v_{i}v_{j}\\
&\overset{(\ref{ee2})}{=}u_{i}v_{j}-u_{j}v_{i}-2u_{i}v_{j}\overset{(\ref{e18})}{=}u_{i}v_{j}-v_{j}u_{i}-2u_{i}v_{j}=0,\\
\{H,v_{i},v_{j}\}&=Hv_{i}v_{j}-Hv_{j}v_{i}+2H\cdot v_{i}v_{j}\overset{(\ref{ee22})}{=}v_{i}v_{j}-v_{j}v_{i}-2v_{i}v_{j}=0.
\end{split}
\end{equation}
Also, note that the function $\{x,y,z\}$ is skew-symmetry in the last two variables, then it is only necessary to have the following equalities:
\begin{equation}\label{lie3}
\begin{split}
\,\,\{E,u_{i},v_{j}\}&=Eu_{i}v_{j}-Ev_{j}u_{i}+2E\cdot u_{i}v_{j}\overset{(\ref{e16})}{=}-v_{i}v_{j}+2(\frac{1}{2}v_{i}v_{j})=0,\\
\,\,\{F,u_{i},v_{j}\}&=Fu_{i}v_{j}-Fv_{j}u_{i}+2F\cdot u_{i}v_{j}\overset{(\ref{e6})}{=}-u_{j}u_{i}+2(-\frac{1}{2}u_{i}u_{j})=0,\\
\,\,\{H,u_{i},v_{j}\}&=Hu_{i}v_{j}-Hu_{j}v_{i}+2H\cdot u_{i}v_{j}\overset{(\ref{e10})}{=}-u_{i}v_{j}-v_{j}u_{i}=0,
\end{split}
\end{equation}
and similarly
\begin{equation}\label{lie4}
\{E,v_{i},u_{j}\}=\{F,v_{i},u_{j}\}=\{H,v_{i},u_{j}\}=0.
\end{equation}
Therefore, by the relations $(\ref{ee24})-(\ref{lie4})$, $J^{2}_{\mathcal{M}}$ is a Lie $L$-module. 

The lemma is proved.
\end{dem}

It follows immediately from Lemma $\ref{l2}$, Proposition $\ref{c1}$ and Lemma $\ref{l3}$ the following result.

\begin{coro}\label{cor1}
    $\mathcal{M}=N_{\mathcal{M}}\oplus J_{\mathcal{M}}$ is a $\mathbf{Z}_{2}$-graded algebra, where $N_{\mathcal{M}}$ is the even part and $J_{\mathcal{M}}$ is the odd part of the $\mathbf{Z}_{2}$-grading of $\mathcal{M}.$
\end{coro}

\subsubsection{Multiplication in $J_{\mathcal{M}}$}

We know that $J_{\mathcal{M}}$ is completely reducible
and is the direct sum of modules isomorphic to the 2-dimensional non-Lie Malcev module, that is, $J_{\mathcal{M}}=\sum_{i} \oplus V_{2i},$ 
where $V_{2i}$ is a non-Lie 2-dimensional Malcev module for the Lie algebra $L.$ Denote by $V(1)$ and $V(2)$ the subspaces of $J_{\mathcal{M}}$ spanned by the elements of
type $u$ and $v$, respectively; then the mappings 
\begin{equation*}
\begin{split}
R_{E}&:V(1)\longrightarrow V(2),\,\,\,\,\,w\mapsto wE,\\
R_{-F}&:V(2)\longrightarrow V(1),\,\,\,\,\,\widetilde{w}\mapsto -\widetilde{w}F\\
\end{split}
\end{equation*}
are mutually inverse and establish isomorphisms between $V(1)$ and $V(2).$ Clearly, $J_{\mathcal{M}}=V(1)\oplus V(2).$ Consider $V=V(1),$ then for any $w\in V$, denote $w(1)=w$ and $w(2)=R_{E}(w).$ Thus $V_{2}(w)=\mathbb{F}\cdot w(1)+\mathbb{F}\cdot w(2)\cong V_{2}.$ So the relations of $(\ref{e5})$ can be written as follows:
\begin{equation}\label{ee35}
    \begin{array}{lclllllll}
w(1)E & = & w(2),\, & w(1)H & = & w(1),\, & w(1)F& = & 0,\vspace{1.5mm}\\
 w(2)E  & = & 0,\, & w(2)H & = & -w(2),\, & w(2)F& = &-w(1),
    \end{array}
\end{equation}
for any $w\in V.$

\begin{pro}\label{p3}
     For any $u, v \in V$, we have
     \begin{center}
         $V_{2}(u)\cdot V_{2}(v)=L\langle u,v	\rangle$
     \end{center}      
where $\langle \cdot,\cdot\rangle:V\times V\longrightarrow U$ is a skew-symmetric bilinear mapping. In particular, $V_{2}(v)^{2}=0$ for any $v\in V.$
\end{pro}

\begin{dem}
    By Corollary \ref{cor1}, $J_{\mathcal{M}}^{2}\subseteq N_{\mathcal{M}}$, then there exists $\alpha_{E}, \alpha_{F}, \alpha_{H}\in U$ such that $u(1)v(1)=E\alpha_{E}+F\alpha_{F}+H\alpha_{H}.$ So by Lemma $\ref{l3}$, $0\overset{(\ref{ee20})}{=}u(1)v(1)F=\frac{1}{2}H\alpha_{E}+F\alpha_{H}$
    which implies $\alpha_{E}=\alpha_{H}=0,$ then $u(1)v(1)=F\alpha_{F}.$
    Denote $\alpha=\alpha_{F}\in U,$ then
    \begin{center}
        $u(1)v(1)=F\alpha.$
    \end{center}
    Furthermore, we have
\begin{equation*}
\begin{split}
u(1)v(2)&\overset{(\ref{ee2})}{=}-u(1)v(1)E=-F\alpha E=\frac{1}{2}H\alpha,\\
u(2)v(1)&\overset{(\ref{ee2})}{=}v(1)u(1)E=-F\alpha E=\frac{1}{2}H\alpha,\\
u(2)v(2)&\overset{(\ref{e16})}{=}-2u(1)v(2)E=-2(\frac{1}{2}H\alpha)E=E\alpha\\
\end{split}
\end{equation*}
which proves that $V_{2}(u)\cdot V_{2}(v)=L\alpha.$ 

Finally, denote $\alpha=\langle u,v\rangle$, then 
\begin{equation*}
V_{2}(u)\cdot V_{2}(v)=L\langle u,v\rangle,\,\,\,\,\,
V_{2}(v)\cdot V_{2}(u)=L\langle v,u \rangle.
\end{equation*}
But from $(\ref{e18})$
\begin{center}
    $u(1)v(2)=-v(1)u(2),$
\end{center}
then $\frac{1}{2}H\langle u,v\rangle =-\frac{1}{2}H\langle v,u\rangle,$ which implies 
    $L(\langle u,v\rangle+\langle v,u\rangle)=0.$
So 
\begin{center}
    $\langle u,v\rangle =-\langle v,u\rangle.$ 
\end{center}
Also
    $F\langle u,v+w\rangle=u(1)(v+w)(1)=u(1)v(1)+u(1)w(1)=F\langle u,v\rangle +F\langle u,w\rangle,$ then
$L(\langle u,v+w\rangle -\langle u,v\rangle -\langle u,w\rangle)=0$ and 
\begin{center}
    $\langle u,v+w\rangle -\langle u,v\rangle -\langle u,w\rangle=0.$
\end{center}
Similarly $\langle u+v,w\rangle -\langle u,v\rangle -\langle u,w\rangle=0,$ which proves that $\langle u,v\rangle$ is a bilinear function of $u$ and $v,$ and the proof is complete.
\end{dem}

\begin{lem}
      For any $u, v, w, t \in V$ the following identities hold:
\begin{equation}\label{ee28}
\,\,\,\,\,\,\,\,\,\,\,\,\,\,\,\,\,\,\,\,\,\,\,\,\,\,\,\,
w\langle u,v\rangle + u\langle v,w\rangle + v\langle w,u\rangle  = 0,
\end{equation}
\begin{equation}\label{ee29}
\langle w,t\rangle\langle u,v\rangle + \langle u,t\rangle\langle v,w\rangle + \langle v,t\rangle\langle w,u\rangle =  0.
\end{equation}
 \end{lem}
\begin{dem} 
Recall that in the proof of Proposition $\ref{p3}$, we obtained the equalities
\begin{equation}\label{ee32}
\begin{split}
u(1)v(1)&=F\langle u,v\rangle,\\
u(1)v(2)=u(2)v(1)&=\frac{1}{2}H\langle u,v\rangle,\\
u(2)v(2)&=E\langle u,v\rangle,\\
\end{split}
\end{equation}
with $\langle u,v\rangle\in U$ and $U$ is a subalgebra of $\Gamma(\mathcal{M}).$ Then $u(1)v(1)\cdot w(2)=F\langle u,v\rangle w(2)=w(1)\langle u,v\rangle$ and
\bes
&u(1)v(1)\cdot w(2)=u(1)v(1)\cdot w(1)E=&\\
&\overset{(\ref{ee1})}{=} u(1)w(1)v(1)E+w(1)v(1)Eu(1)+v(1)Eu(1)w(1)+Eu(1)w(1)v(1)=&\\
&\overset{(\ref{ee35}),(\ref{ee32})}{=}F\langle u,w\rangle v(1)E+F\langle w,v\rangle Eu(1)+v(2)u(1)w(1)-u(2)w(1)v(1)=&\\
&=-\frac{1}{2} H\langle w,v\rangle u(1)+\frac{1}{2}H\langle v,u\rangle w(1)-\frac{1}{2}H\langle u,w\rangle v(1)=&\\
&=\frac{1}{2}u(1)\langle w,v\rangle -\frac{1}{2}w(1)\langle v,u\rangle +\frac{1}{2}v(1)\langle u,w\rangle,
\ees
thus $w(1)\langle u,v\rangle =\frac{1}{2} u(1)\langle w,v\rangle - \frac{1}{2} w(1)\langle v,u\rangle+\frac{1}{2}v(1)\langle u,w\rangle$ and 
\begin{equation}\label{super}
    w(1)\langle u,v\rangle = u(1)\langle w,v\rangle + v(1)\langle u,w\rangle,
\end{equation}
which proves $(\ref{ee28})$. Multiplying $(\ref{ee28})$ by the element $t\in V,$ we get $(\ref{ee29}).$ The lemma is proved.
\end{dem}

Similarly, $-u(2)v(2)\cdot w(1)=-E\langle u,v\rangle w(1)= w(2)\langle u,v\rangle$ and
\bes
        &-u(2)v(2)\cdot w(1)=u(2)v(2)\cdot w(2)F=&\\
        &\overset{(\ref{ee1})}{=} u(2)w(2)v(2)F+w(2)v(2)Fu(2)+v(2)Fu(2)w(2)+Fu(2)w(2)v(2)=&\\
        &\overset{(\ref{ee35}),(\ref{ee32})}{=}E\langle u,w\rangle v(2)F+E\langle w,v\rangle Fu(2)-v(1)u(2)w(2)+u(1)w(2)v(2)=&\\
        &=\frac{1}{2} H\langle w,v\rangle u(2)-\frac{1}{2}H\langle v,u\rangle w(2)+\frac{1}{2}H\langle u,w\rangle v(2)=&\\
        &=\frac{1}{2} u(2)\langle w,v\rangle-\frac{1}{2}w(2)\langle v,u \rangle +\frac{1}{2}v(2)\langle u,w \rangle,
\ees
then $w(2)\langle u,v\rangle =\frac{1}{2} u(2)\langle w,v\rangle-\frac{1}{2}w(2)\langle v,u \rangle +\frac{1}{2}v(2)\langle u,w \rangle$ and 
\begin{equation}\label{ee30}
    w(2)\langle u,v\rangle = u(2)\langle w,v\rangle +v(2)\langle u,w \rangle.
\end{equation}
Therefore, from $(\ref{super})$ and $(\ref{ee30})$, we obtain
\begin{center}
    $V_{2}(w)\langle u,v\rangle \subseteq V_{2}(u)\langle w,v\rangle +V_{2}(v)\langle u,w \rangle.$
\end{center}

The following results were proved in \cite{LSS1}. 

\begin{coro}
Let $\left\{v_i \mid i \in I\right\}$ be a basis of the space $V$ and let $u_{i j}=\left\langle v_i, v_j\right\rangle \in U$. Then the elements $u_{i j}$ satisfy the Plücker relations
\begin{equation}\label{ee31}
    u_{i j}=-u_{j i}, \quad u_{i j} u_{k l}+u_{i k} u_{l j}+u_{i l} u_{j k}=0
\end{equation}
\end{coro}

\begin{lem}
Consider the algebra of polynomials $\mathbb{F}\left[x_1, \ldots, x_n ; y_1, \ldots, y_n\right]$, and let $\alpha_{i j}=\operatorname{det}\left[\begin{array}{ll}x_i & y_i \\ x_j & y_j\end{array}\right] \in \mathbb{F}\left[x_1, \ldots, x_n ; y_1, \ldots, y_n\right]$. Then the elements $\alpha_{i j}=-\alpha_{j i}$ satisfy relations $(\ref
{ee31})$. Moreover, the algebra $\mathbb{F}\left[\alpha_{i j} \mid 1 \leq i<j \leq n\right]$ is a free algebra modulo relations $(\ref{ee31})$.
\end{lem}

In the Lie algebra $\mathfrak{sl}_{2}(\mathbb{F}),$ consider a basis formed by matrices 
\begin{equation*}
E=
\begin{pmatrix}
0 & -\frac{1}{2} \\
0 & 0 \\
\end{pmatrix},\,\,\,\,
H=
\begin{pmatrix}
-\frac{1}{2} & 0 \\
0 & \frac{1}{2} \\
\end{pmatrix}\,\,\
\textup{and}\,\,\
F=
\begin{pmatrix}
0 & 0 \\
\frac{1}{2} & 0 \\
\end{pmatrix}.
\end{equation*}
that clearly satisfy the multiplication table of $\mathfrak{sl}_{2}(\mathbb{F})$ (see the subsection \ref{subsection 2.1}).
Also, recall that by Proposition \ref{c1}, $N_{\mathcal{M}}\cong \mathfrak{sl}_{2}(U),$ hence $\mathcal{M}=\mathfrak{sl}_{2}(U)\oplus V(1)\oplus V(2).$

\begin{pro}\label{p4}
Let $X, Y \in \mathcal{M},$ $X=X_{L}+x(1)+y(2),$ $Y=Y_{L}+z(1)+t(2)$, where $X_{L}=\frac{1}{2}\left(\begin{array}{ll}-c & -a \\ \,\,\,\,b & \,\,\,\,c\end{array}\right),$ $Y_{L}=\frac{1}{2}\left(\begin{array}{ll}-f & -d \\ \,\,\,\,e & \,\,\,\,f\end{array}\right),$ $a, b, c, d, e, f\in U,$ $x, y, z, t \in V$. Then the product $X Y$ is given by
\begin{equation*}
    \begin{split}
X Y&=X_{L} Y_{L}+\dfrac{1}{2}\left(\begin{array}{cc}
-\frac{1}{2}(\langle x, t\rangle+\langle y, z\rangle)  & -\langle y, t\rangle \\
\langle x,z\rangle & \frac{1}{2}(\langle x, t\rangle+\langle y, z\rangle)
\end{array}\right)\\
&\,\,\,\,\,\,+(fx-ey-cz+bt)(1)+(dx-fy-az+ct)(2),
\end{split}
\end{equation*}
where $X_{L} Y_{L}$ is the Lie bracket of $X_{L}$ and $Y_{L}.$
\end{pro}
\begin{dem} 
Directly from the relations $(\ref{ee35})$ and $(\ref{ee32}).$ 
\end{dem} 

We can make the formula defining the product in $\mathcal{M}$ more transparent by using the following notation: 
\begin{equation*}
(u, v)=u(1)+v(2).
\end{equation*}
So, we have $X=X_{L}+(x,y)$ and $Y=Y_{L}+(z,t)$, then 
\begin{equation}\label{ee33}
\begin{split}
X Y&=X_{L} Y_{L}+\dfrac{1}{2}\left(\begin{array}{cc}
-\frac{1}{2}(\langle x, t\rangle+\langle y, z\rangle)  & -\langle y,t\rangle \\
\langle x,z\rangle & \frac{1}{2}(\langle x, t\rangle+\langle y, z\rangle)
\end{array}\right) + (x,y)Y_{L}+X_{L}(z,t).
\end{split}
\end{equation}

In the next subsection we will prove that Proposition \ref{p4} in fact describes all Malcev extensions of the Lie algebra $\mathfrak{sl}_{2}(\mathbb{F})$ in the variety $\mathcal{H}.$ 

\subsubsection{The main theorem}\label{subsection}

Let $\mathcal{B}$ be a unital associative commutative algebra over a field of characteristic $\neq 2,3$ and let $V$ be a commutative $\mathcal{B}$-bimodule. Assume that there exists a $\mathcal{B}$-bilinear skew-symmetric mapping $\langle\cdot , \cdot\rangle: V \times V \rightarrow \mathcal{B}$ such that $\langle V, V\rangle \subseteq \mathcal{B}$ and formula $(\ref{ee28})$ holds for any $u, v, w \in V$.

Let $\mathcal{M}=\mathfrak{sl}_{2}(\mathcal{B}) \oplus V^2$, where $\mathfrak{sl}_{2}(\mathcal{B})$ is the vector space of all matrices “coordinated” by $\mathcal{B}$ having trace zero, $V^2=\{(u, v) \mid u, v \in V\} \cong V \oplus V$ and
\begin{center}
    $\mathcal{B}\,\,\textup{a subalgebra of}\,\, \Gamma(\mathcal{M}).$
\end{center}
Let $X, Y \in \mathcal{M}$, $X=X_{L}+(x, y), Y=Y_{L}+(z, t)$, where $X_{L}, Y_{L} \in \mathfrak{sl}_{2}(\mathcal{B})$ and $(x, y),(z, t) \in V^2$. Define a product in $\mathcal{M}$ by formula $(\ref{ee33})$:
\begin{equation*}
\begin{split}
X Y&=X_{L} Y_{L}+\dfrac{1}{2}\left(\begin{array}{cc}
-\frac{1}{2}(\langle x, t\rangle+\langle y, z\rangle)  & -\langle y,t\rangle \\
\langle x,z\rangle & \frac{1}{2}(\langle x, t\rangle+\langle y, z\rangle)
\end{array}\right) + (x,y)Y_{L}+X_{L}(z,t).
\end{split}
\end{equation*}
where $X_{L} Y_{L}$ is the Lie bracket of $X_{L}$ and $Y_{L}$.
\begin{theo}\label{principal}
The algebra $\mathcal{M}$ with the product defined above is a Malcev algebra in $\mathcal{H}$ containing $L=\mathfrak{sl}_{2}.$ 
Conversely, every Malcev algebra in $\mathcal{H}$ that contains a subalgebra $L\cong\mathfrak{sl}_{2}$ with $mL\neq 0$ for any $0\neq m\in\mathcal{M}$ has this form.
\end{theo}
\begin{dem}
The second part of the theorem follows from Proposition 
\ref{p4} with $\mathcal{B}=U$. For the first part, we begin proving that $\mathcal{M}$ is a Malcev algebra.

Let $B$ be a unital associative algebra and let $W$ be a right $B$-module such that $W[B,B]=0$, then in this case $w\cdot b=b\cdot w$ for all $w\in W$, $b\in B.$ Assume that there exist a $B$-bilinear skew-symmetric mapping 
$\langle\cdot , \cdot\rangle: W \times W \rightarrow B$ such that $\langle W, W\rangle\subseteq Z(B)$ and for any $u,v,w\in W$
\begin{equation*}
    \langle u,v\rangle w+\langle v,w\rangle u+ \langle w,u\rangle v=0.
\end{equation*}

Let $A=M_2(B)\oplus W^2,$ where $W^2=\{(u,v)|\,u,v\in W\}\cong W\oplus W.$ Let $X,Y\in A,$ $X=X_a+(x,y),$ $Y=Y_a+(z,t),$ where $X_a,Y_a\in M_2(B)$ and $(x,y), (z,t)\in W^2.$
Define a product in $A$ by formula:
\begin{equation*}
    XY=X_aY_a+\frac{1}{4}\left(\begin{array}{cc}
-\langle x,t\rangle  & -\langle y,t\rangle \\
\langle x,z\rangle & \langle y,z\rangle
\end{array}\right)+(z,t)(X_a)^{\star}+(x,y)Y_a,
\end{equation*}
where $A \mapsto A^*$ is the symplectic involution in $M_2(B)$:
\begin{center}
$\left(\begin{array}{cc}
a & b \\
c & d
\end{array}\right)^*
\mapsto
\left(\begin{array}{cc}
d & -b \\
-c & a
\end{array}\right)$.
\end{center}
Using the proof of Theorem 5.1 of \cite{LSS1}, it is prove that $A$ is a unital alternative algebra. From this we know that $A^{-}=(A,[\,,\,])$ is a Malcev algebra, where $[\,,\,]$ is the Lie bracket. The product in $A^{(-)}$ is: 
\bes
&[X,Y]=XY-YX=&
\ees
\bes
&=X_aY_a+\frac{1}{4}\left(\begin{array}{cc}
-\langle x,t\rangle  & -\langle y,t\rangle \\
\langle x,z\rangle & \langle y,z\rangle
\end{array}\right)+(z,t)(X_a)^{\star}+(x,y)Y_a-&\\
&-Y_aX_a-\frac{1}{4}\left(\begin{array}{cc}
-\langle z,y\rangle  & -\langle t,y\rangle \\
\langle z,x\rangle & \langle t,x\rangle
\end{array}\right)-(x,y)(Y_a)^{\star}-(z,t)X_a=&\\
&=[X_a,Y_a]+\dfrac{1}{2}\left(\begin{array}{cc}
-\frac{1}{2}(\langle x, t\rangle+\langle y, z\rangle)  & -\langle y,t\rangle \\
\langle x,z\rangle & \frac{1}{2}(\langle x, t\rangle+\langle y, z\rangle)
\end{array}\right)+[(x,y),Y_a]+[X_a,(z,t)].
\ees
Take $\mathcal{B}=Z(B)$ and $R=M_2(\mathcal{B})\oplus W^2.$ Thus we see that $\mathcal{M}\cong R^{-} \leq A^{-}$, which shows that $\mathcal{M}$ isomorphic to a subalgebra of  $A^{-}.$ Then $\mathcal{M}$ is a Malcev algebra.

Finally, we will prove that $\mathcal{M}$ belongs to the variety $\mathcal{H}$. For this it will be enough to verify that $\alpha(\mathcal{M},\mathcal{M},\mathcal{M})$ is contained in $\Gamma(\mathcal{M})$ (see \cite{F1}). It is clear that $\alpha(\mathfrak{sl}_2(\mathcal{B}),\mathfrak{sl}_2(\mathcal{B}),\mathfrak{sl}_2(\mathcal{B}))$ is contained in $\Gamma(\mathcal{M}).$ Now as $\alpha$ is skew-symmetric, we consider only following cases:
\begin{center}
    \textbf{The cases $\alpha((u,v),E,F)$, $\alpha((u,v),E,H)$, $\alpha((u,v),F,H)$}
\end{center}
We begin with
\bes
&E\alpha((u,v),E,F)=p(E,(u,v),E,F)=&\\
&=-\{EF,E,(u,v)\}-\{(u,v)F,E,E\}+\{EF,(u,v),E\}=&\\
&=-\frac{1}{2}\{H,E,(u,v)\}+\{v(1),E,E\}+\frac{1}{2}\{H,(u,v),E\}=&\\
&=-\frac{1}{2}(HE(u,v)-H(u,v)E+ 2H\cdot E(u,v))+&\\
&+\frac{1}{2}(H(u,v)E-HE(u,v)+2H\cdot (u,v)E)=&\\
&=-\frac{1}{2}(-E(u,v)-(-u(1)+v(2))E-2H u(2))+\frac{1}{2}((-u(1)+v(2))E+&\\
&+E(u,v)+2Hu(2))=&\\
&=E(u,v)+(-u(1)+v(2))E+2H u(2)=-u(2)-u(2)+2u(2)=0.&
\ees
Also, 
\bes
&F\alpha((u,v),E,F)=p(F,(u,v),E,F)=&\\
&=-\{EF,F,(u,v)\}-\{(u,v)F,E,F\}+\{FF,(u,v),E\}&\\
&=-\frac{1}{2}\{H,F,(u,v)\}+\{v(1),E,F\}=&\\
&=-\frac{1}{2}(HF(u,v)-H(u,v)F+ 2H\cdot F(u,v))+v(1)EF-v(1)FE+2v(1)\cdot EF=&\\
&=-\frac{1}{2}(F(u,v)-(-u(1)+v(2))F+2Hv(1))+v(2)F+v(1)H=&\\
&=-\frac{1}{2}(v(1)+v(1)-2v(1))-v(1)+v(1)=0.&
\ees
In the same way,
\bes
&H\alpha((u,v),E,F)=p(H,(u,v),E,F)=-\{EF,H,(u,v)\}-\{(u,v)F,E,H\}+&\\
&+\{HF,(u,v),E\}=-\frac{1}{2}\{H,H,(u,v)\}+\{v(1),E,H\}+\{F,(u,v),E\}=&\\
&=-\frac{1}{2}(HH(u,v)-H(u,v)H+ 2H\cdot H(u,v))+v(1)EH- v(1)HE+2v(1)\cdot EH +&\\
&+F(u,v)E-FE(u,v)+2F\cdot (u,v)E=-\frac{1}{2}(-(-u(1)+v(2))H+2H(-u(1)+v(2)))+&\\
&+v(2)H-v(1)E+2v(1)E+v(1)E+\frac{1}{2}H(u,v)+2Fu(2)=0.&
\ees
Besides,
\bes
&x(1)\alpha((u,v),E,F)=p(x(1),(u,v),E,F)=-\{EF,x(1),(u,v)\}-\{(u,v)F,E,x(1)\}+&\\
&+\{x(1)F,(u,v),E\}=-\frac{1}{2}\{H,x(1),(u,v)\}-\{v(1),E,x(1)\}=&\\
&=-\frac{1}{2}(Hx(1)(u,v)-H(u,v)x(1)+2H\cdot x(1)(u,v))-&\\
&-(v(1)Ex(1)-v(1)x(1)E+2v(1)\cdot Ex(1))=
&\\
&=-\frac{1}{2}(-x(1)(u,v)-(-u(1)+v(2))x(1)+2H(F\langle x,u\rangle+\frac{1}{2}H\langle x,v\rangle))-&\\
&\,-(v(2)x(1)-F\langle v,x\rangle E-2v(1)x(2))=\frac{1}{2}(F\langle x,u\rangle+\frac{1}{2}H\langle x,v\rangle-F\langle u,x\rangle+\frac{1}{2}H\langle v,x\rangle -&\\
&-2F\langle x,u\rangle)-\frac{1}{2}H\langle v,x\rangle-\frac{1}{2}H\langle v,x\rangle+H\langle v,x\rangle=0&
\ees
and
\bes
&y(2)\alpha((u,v),E,F)=p(y(2),(u,v),E,F)=-\{EF,y(2),(u,v)\}-\{(u,v)F,E,y(2)\}+&\\
&+\{y(2)F,(u,v),E\}=-\frac{1}{2}\{H,y(2),(u,v)\}+\{v(1),E,y(2)\} - \{y(1),(u,v),E\}=&\\
&=-\frac{1}{2}(Hy(2)(u,v)-H(u,v)y(2)+2H\cdot y(2)(u,v))+&\\
&+v(1)Ey(2)-v(1)y(2)E+2v(1)\cdot Ey(2)-&\\
&-y(1)(u,v)E+y(1)E(u,v)-2y(1)\cdot (u,v)E=&\\
&=-\frac{1}{2}(y(2)(u,v)-(-u(1)+v(2))y(2)+2H(\frac{1}{2}H\langle y,u\rangle+E\langle y,v\rangle))+&\\
&+v(2)y(2)-\frac{1}{2}H\langle v,y\rangle E-(F\langle y,u\rangle+\frac{1}{2}H\langle y,v\rangle)E+y(2)(u,v)-2y(1)u(2)=&\\
&=-\frac{1}{2}(y(2)u(1)+y(2)v(2)+u(1)y(2)-v(2)y(2)-2E\langle y,v\rangle)+&\\
&+v(2)y(2)+\frac{1}{2}E\langle v,y\rangle+\frac{1}{2}H\langle y,u\rangle+\frac{1}{2}E\langle y,v\rangle+y(2)u(1)+y(2)v(2)-H\langle y,u\rangle=&\\
&=\frac{1}{2}E\langle v,y\rangle+\frac{1}{2}H\langle y,u\rangle+\frac{1}{2}E\langle y,v\rangle+\frac{1}{2}H\langle y,u\rangle-H\langle y,u\rangle =0. &
\ees
Similarly 
    $\alpha((u,v),E,H)|_\mathcal{M}=0$ and $\alpha((u,v),F,H)|_\mathcal{M}=0,$
then obviously $\alpha((u,v),E,F),$ $\alpha((u,v),E,H)$ and $\alpha((u,v),F,H)$ belong to $\Gamma(\mathcal{M}).$
\begin{center}
    \textbf{The cases $\alpha((u,v),(z,t),E)$, $\alpha((u,v),(z,t),F)$, $\alpha((u,v),(z,t),H)$}
\end{center}
Now we begin with
\bes
&E\alpha((u,v),(z,t),E)=p(E,(u,v),(z,t),E)=&\\
&=-\{(z,t)E,E,(u,v)\}-\{(u,v)E,(z,t),E\}+\{EE,(u,v),(z,t)\}=&\\
&=-\{z(2),E,(u,v)\}-\{u(2),(z,t),E\}=&\\
&=-z(2)E(u,v)+z(2)(u,v)E-2z(2)\cdot E(u,v)-&\\
&-u(2)(z,t)E+u(2)E(z,t)-2u(2)\cdot (z,t)E=&\\
&=(\frac{1}{2}H\langle z,u\rangle+E\langle z,v\rangle)E+2z(2)u(2)-(\frac{1}{2}H\langle u,z\rangle+E\langle u,t\rangle)E-2u(2)z(2)=&\\
&=-\frac{1}{2}E\langle z,u\rangle+2E\langle z,u\rangle+\frac{1}{2}E\langle u,z\rangle-2E\langle u,z\rangle=3E\langle z,u\rangle.&
\ees
Furthermore,
\bes
&F\alpha((u,v),(z,t),E)=p(F,(u,v),(z,t),E)=&\\
&=-\{(z,t)E,F,(u,v)\}-\{(u,v)E,(z,t),F\}+\{FE,(u,v),(z,t)\}=\\
&=-\{z(2),F,(u,v)\} -\{u(2),(z,t),F\} -\frac{1}{2}\{H,(u,v),(z,t)\}=&\\
&=-z(2)F(u,v)+z(2)(u,v)F-2z(2)\cdot F(u,v)-&\\
&- u(2)(z,t)F+u(2)F(z,t)-2u(2)\cdot (z,t)F-&\\
&-\frac{1}{2}(H(u,v)(z,t)-H(z,t)(u,v)+2H\cdot (u,v)(z,t))=&\\
&=z(1)(u,v)+(\frac{1}{2}H\langle z,u\rangle + E\langle z,v\rangle)F-2z(2)v(1)-&\\
&-(\frac{1}{2}H\langle u,z\rangle+E\langle u,t\rangle)F-u(1)(z,t)+2u(2)t(1)-&\\
&-\frac{1}{2}((-u(1)+v(2))(z,t)-(-z(1)+t(2))(u,v)+&\\
&+2H(F\langle u,z\rangle+\frac{H}{2}\langle u, t\rangle+\frac{H}{2}\langle v, z\rangle+E\langle v,t\rangle))=&\\
&=F\langle z,u\rangle+\frac{1}{2}H\langle z,v\rangle+\frac{1}{2}F\langle z,u\rangle + \frac{1}{2}H\langle z,v\rangle - H\langle z,v\rangle-&\\
&-\frac{1}{2}F\langle u,z\rangle-\frac{1}{2}H\langle u,t\rangle - F\langle u,z\rangle - \frac{1}{2}H \langle u,t\rangle + H\langle u,t\rangle-&\\
&-\frac{1}{2}(-F\langle u,z\rangle-\frac{1}{2}H\langle u,t\rangle+\frac{1}{2}H\langle v,z\rangle + E\langle v,t\rangle+&\\
&+F\langle z,u\rangle+\frac{1}{2}H\langle z,v\rangle-\frac{1}{2}H\langle t,u\rangle-E\langle t,v\rangle +&\\
&+2F\langle u,z\rangle - 2E\langle v,t\rangle)=3F\langle z,u\rangle.
\ees
Also,
\bes
&H\alpha((u,v),(z,t),E)=p(H,(u,v),(z,t),E)=&\\
&=-\{(z,t)E,H,(u,v)\}-\{(u,v)E,(z,t),H\}+\{HE,(u,v),(z,t)\}=&\\
&=-\{z(2),H,(u,v)\}-\{u(2),(z,t),H\}-\{E,(u,v),(z,t)\}=&\\
&=-(z(2)H(u,v)-z(2)(u,v)H+2z(2)\cdot H(u,v))-&\\
&-(u(2)(z,t)H-u(2)H(z,t)+2u(2)\cdot (z,t)H)-&\\
&-(E(u,v)(z,t)-E(z,t)(u,v)+2E\cdot (u,v)(z,t))=&\\
&=z(2)(u,v)+(\frac{1}{2}H\langle z,u\rangle+E\langle z,v\rangle)H-2z(2)(-u(1)+v(2))-&\\
&-(\frac{1}{2}H\langle u,z\rangle + E\langle u,t\rangle)H-u(2)(z,t)-2u(2)(z(1)-t(2))+&\\
&+u(2)(z,t)-z(2)(u,v)-2E(F\langle u,z\rangle+\frac{H}{2}\langle u, t\rangle+\frac{H}{2}\langle v, z\rangle+E\langle v,t\rangle)=&\\
&=\frac{1}{2}H\langle z,u\rangle+E\langle z,v\rangle + E\langle z,v\rangle + H\langle z,u\rangle - 2E\langle z,v\rangle-&\\
&-E\langle u,t\rangle-\frac{1}{2}H\langle u,z\rangle-E\langle u,t\rangle-H\langle u,z\rangle+2E\langle u,t\rangle+&\\
&+\frac{1}{2}H\langle u,z\rangle + E\langle u,t\rangle-\frac{1}{2}H\langle z,u\rangle - E\langle z,v\rangle -&\\
&-H\langle u,z\rangle - E\langle u,t\rangle - E\langle v,z\rangle = 3H\langle z,u\rangle.&
\ees
In addition, 
\bes
&x(1)\alpha((u,v),(z,t),E)=p(x(1),(u,v),(z,t),E)=&\\
&=-\{(z,t)E,x(1),(u,v)\}-\{(u,v)E,(z,t),x(1)\}+\{x(1)E,(u,v),(z,t)\}=&\\
&=-\{z(2),x(1),(u,v)\}-\{u(2),(z,t),x(1)\}+\{x(2),(u,v),(z,t)\}=&\\
&=-(z(2)x(1)(u,v)-z(2)(u,v)x(1)+2z(2)\cdot x(1)(u,v))-&\\
&-(u(2)(z,t)x(1)-u(2)x(1)(z,t)+2u(2)\cdot (z,t)x(1))+&\\
&+x(2)(u,v)(z,t)-x(2)(z,t)(u,v)+2x(2)\cdot (u,v)(z,t)&\\
&=-\frac{1}{2}H\langle z,x\rangle (u,v)+(\frac{1}{2}H\langle z,u\rangle + E\langle z,v\rangle)x(1)-2z(2)(F\langle x,u\rangle+\frac{1}{2}H\langle x,v\rangle)-&\\
&-(\frac{1}{2}H\langle u,z\rangle + E\langle u,t\rangle)x(1) + \frac{1}{2}H\langle u,x\rangle (z,t) - 2u(2)(F\langle z,x\rangle+\frac{1}{2}H\langle t,x\rangle)+&\\
&+(\frac{1}{2}H\langle x,u\rangle + E\langle x,v\rangle)(z,t)-(\frac{1}{2}H\langle x,z\rangle+E\langle x,t\rangle)(u,v)+&\\
&+2x(2)(F\langle u,z\rangle+\frac{H}{2}\langle u, t\rangle+\frac{H}{2}\langle v, z\rangle+E\langle v,t\rangle)=&\\
&=\frac{1}{2}u(1)\langle z,x\rangle - \frac{1}{2}v(2)\langle z,x\rangle -\frac{1}{2}x(1)\langle z,u\rangle-x(2)\langle z,v\rangle + 2z(1)\langle x,u\rangle+z(2)\langle x,v\rangle+&\\
&+\frac{1}{2}x(1)\langle u,z\rangle+x(2)\langle u,t\rangle-\frac{1}{2}z(1)\langle u,x\rangle+\frac{1}{2}t(2)\langle u,x\rangle + 2u(1)\langle z,x\rangle+u(2)\langle t,x\rangle-&\\
&-\frac{1}{2}z(1)\langle x,u\rangle+\frac{1}{2}t(2)\langle x,u\rangle-z(2)\langle x,v\rangle +&\\
&+\frac{1}{2}u(1) \langle x,z\rangle -\frac{1}{2}v(2)\langle x,z\rangle+u(2)\langle x,t\rangle-&\\
&-2x(1)\langle u,z\rangle-x(2)\langle u,t\rangle - x(2)\langle v,z\rangle=&\\
&=2(u(1)\langle z,x\rangle+z(1)\langle x,u\rangle)-x(1)\langle u,z\rangle\overset{(\ref{ee28})}{=}3x(1)\langle z,u\rangle&
\ees
and
\bes
&y(2)\alpha((u,v),(z,t),E)=p(y(2),(u,v),(z,t),E)=&\\
&=-\{(z,t)E,y(2),(u,v)\}-\{(u,v)E,(z,t),y(2)\}+\{y(2)E,(u,v),(z,t)\}=&\\
&=-\{z(2),y(2),(u,v)\}-\{u(2),(z,t),y(2)\}=&\\
&=-(z(2)y(2)(u,v)-z(2)(u,v)y(2)+2z(2)\cdot y(2)(u,v))-&\\
&-(u(2)(z,t)y(2)-u(2)y(2)(z,t)+2u(2)\cdot (z,t)y(2))=&\\
&=-E\langle z,y\rangle (u,v)+(\frac{1}{2}H\langle z,u\rangle+E\langle z,v\rangle)y(2)-2z(2)(\frac{1}{2}H\langle y,u\rangle+E\langle y,v\rangle)-&\\
&-(\frac{1}{2}H\langle u,z\rangle + E\langle u,t\rangle)y(2)+E\langle u,y\rangle (z,t)-2u(2)(\frac{1}{2}H\langle z,y\rangle + E\langle t,y\rangle)=&\\
&=u(2)\langle z,y\rangle+\frac{1}{2}y(2)\langle z,u\rangle+z(2)\langle y,u\rangle-\frac{1}{2}y(2)\langle u,z\rangle - z(2)\langle u,y\rangle + u(2)\langle z,y\rangle =&\\
&=2(u(2)\langle z,y\rangle +z(2)\langle y,u\rangle )+y(2)\langle z,u\rangle \overset{(\ref{ee28})}{=}3y(2)\langle z,u\rangle.& 
\ees
Similarly,
\begin{equation*}
    \begin{split}
        E\alpha((u,v),(z,t),F)&=3E\langle t,v\rangle,\\
        F\alpha((u,v),(z,t),F)&=3F\langle t,v\rangle,\\
        H\alpha((u,v),(z,t),F)&=3H\langle t,v\rangle,\\
        x(1)\alpha((u,v),(z,t),F)&=3x(1)\langle t,v\rangle,\\
        y(2)\alpha((u,v),(z,t),F)&=3y(2)\langle t,v\rangle.
    \end{split}
\end{equation*}
Moreover,
\bes
&E\alpha((u,v),(z,t),H)=p(E,(u,v),(z,t),H)=&\\
&=-\{(z,t)H,E,(u,v)\}-\{(u,v)H,(z,t),E\}+\{EH,(u,v),(z,t)\}=&\\
&\overset{(\ref{lie1})-(\ref{lie4})}{=}-\{z(1)-t(2),E,(u,v)\}-\{u(1)-v(2),(z,t),E\}=&\\
&=-\{z(1),E,(u,v)\}+\{t(2),E,(u,v)\}-\{u(1),(z,t),E\}+\{v(2),(z,t),E\}=&\\
&=-z(1)E(u,v)+z(1)(u,v)E-2z(1)\cdot E(u,v)+&\\
&+t(2)E(u,v)-t(2)(u,v)E+2t(2)\cdot E(u,v)-&\\
&-u(1)(z,t)E+u(1)E(z,t)-2u(1)\cdot (z,t)E+&\\
&+v(2)(z,t)E-v(2)E(z,t)+2v(2)\cdot (z,t)E=&\\
&=-z(2)(u,v)+(F\langle z,u\rangle+\frac{1}{2}H\langle z,v\rangle)E+2z(1)u(2)-&\\
&-(\frac{1}{2}H\langle t,u\rangle+E\langle t,v\rangle)E-2t(2)u(2)-&\\
&-(F\langle u,z\rangle+\frac{1}{2}H\langle u,t\rangle)E+u(2)(z,t)-2u(1)z(2)+&\\
&+(\frac{1}{2}H\langle v,z\rangle+E\langle v,t\rangle)E+2v(2)z(2)=&\\
&=-\frac{1}{2}H\langle z,u\rangle-E\langle z,v\rangle-\frac{1}{2}H\langle z,u\rangle -\frac{1}{2}E\langle z,v\rangle +H\langle z,u\rangle+
+\frac{1}{2}E\langle t,u\rangle-2E\langle t,u\rangle+&\\
&+\frac{1}{2}H\langle u,z\rangle +\frac{1}{2}E\langle u,t\rangle +\frac{1}{2}H\langle u,z\rangle +E\langle u,t\rangle -H\langle u,z\rangle
-\frac{1}{2}E\langle v,z\rangle+2E\langle v,z\rangle=&\\
&=3E(\langle u,t\rangle+\langle v,z\rangle).&
\ees
Further,
\bes
&F\alpha((u,v),(z,t),H)=p(F,(u,v),(z,t),H)=&\\
&=-\{(z,t)H,F,(u,v)\}-\{(u,v)H,(z,t),F\}+\{FH,(u,v),(z,t)\}=&\\
&\overset{(\ref{lie1})-(\ref{lie4})}{=}-\{z(1)-t(2),F,(u,v)\}-\{u(1)-v(2),(z,t),F\}=&\\
&=-\{z(1),F,(u,v)\}+\{t(2),F,(u,v)\}-\{u(1),(z,t),F\}+\{v(2),(z,t),F\}=&\\
&=-z(1)F(u,v)+z(1)(u,v)F-2z(1)\cdot F(u,v)+&\\
&+t(2)F(u,v)-t(2)(u,v)F+2t(2)\cdot F(u,v)-&\\
&-u(1)(z,t)F+u(1)F(z,t)-2u(1)\cdot (z,t)F+&\\
&+v(2)(z,t)F-v(2)F(z,t)+2v(2)\cdot (z,t)F=&\\
&=(F\langle z,u\rangle+\frac{1}{2}H\langle z,v\rangle)F-2z(1)v(1)-&\\
&-t(1)(u,v)-(\frac{1}{2}H\langle t,u\rangle+E\langle t,v\rangle)F+2t(2)v(1)-&\\
&-(F\langle u,z\rangle+\frac{1}{2}H\langle u,t\rangle)F+2u(1)t(1)+&\\
&+(\frac{1}{2}H\langle v,z\rangle+E\langle v,t\rangle)F+v(1)(z,t)-2v(2)t(1)=&\\
&=\frac{1}{2}F\langle z,v\rangle-2F\langle z,v\rangle-F\langle t,u\rangle-\frac{1}{2}H\langle t,v\rangle-\frac{1}{2}F\langle t,u\rangle-\frac{1}{2}H\langle t,v\rangle+H\langle t,v\rangle -&\\
&-\frac{1}{2}F\langle u,t\rangle+2F\langle u,t\rangle+\frac{1}{2}F\langle v,z\rangle+\frac{1}{2}H\langle v,t\rangle+F\langle v,z\rangle+\frac{1}{2}H\langle v,t\rangle-H\langle v,t\rangle=&\\
&=3F(\langle u,t\rangle+\langle v,z\rangle).&
\ees
The same way,
\bes
&H\alpha((u,v),(z,t),H)=p(H,(u,v),(z,t),H)=&\\
&=-\{(z,t)H,H,(u,v)\}-\{(u,v)H,(z,t),H\}+\{HH,(u,v),(z,t)\}=&\\
&=-\{z(1)-t(2),H,(u,v)\}-\{u(1)-v(2),(z,t),H\}=&\\
&=-\{z(1),H,(u,v)\}+\{t(2),H,(u,v)\}-\{u(1),(z,t),H\}+\{v(2),(z,t),H\}=&\\
&=-z(1)H(u,v)+z(1)(u,v)H-2z(1)\cdot H(u,v)+&\\
&+t(2)H(u,v)-t(2)(u,v)H+2t(2)\cdot H(u,v)-&\\
&-u(1)(z,t)H+u(1)H(z,t)-2u(1)\cdot (z,t)H+&\\
&+v(2)(z,t)H-v(2)H(z,t)+2v(2)\cdot (z,t)H=&\\
&=-z(1)(u,v)+(F\langle z,u\rangle+\frac{1}{2}H\langle z,v\rangle)H-2z(1)(-u(1)+v(2))-&\\
&-t(2)(u,v)-(\frac{1}{2}H\langle t,u\rangle+E\langle t,v\rangle)H+2t(2)(-u(1)+v(2))-&\\
&-(F\langle u,z\rangle+\frac{1}{2}H\langle u,t\rangle)H+u(1)(z,t)-2u(1)(z(1)-t(2))+&\\
&+(\frac{1}{2}H\langle v,z\rangle+E\langle v,t\rangle)H+v(2)(z,t)+2v(2)(z(1)-t(2))=&\\
&=-F\langle z,u\rangle-\frac{1}{2}H\langle z,v\rangle-F\langle z,u\rangle+2F\langle z,u\rangle-H\langle z,v\rangle-&\\
&-\frac{1}{2}H\langle t,u\rangle-E\langle t,v\rangle-E\langle t,v\rangle-H\langle t,u\rangle+2E\langle t,v\rangle+&\\
&+F\langle u,z\rangle+F\langle u,z\rangle+\frac{1}{2}H\langle u,t\rangle-2F\langle u,z\rangle+H\langle u,t\rangle +&\\
&+E\langle v,t\rangle+\frac{1}{2}H\langle v,z\rangle+E\langle v,t\rangle+H\langle v,z\rangle-2E\langle v,t\rangle=&\\
&=3H(\langle u,t\rangle+\langle v,z\rangle).&
\ees
Following the same process,
\bes
&x(1)\alpha((u,v),(z,t),H)=p(x(1),(u,v),(z,t),H)=&\\
&=-\{(z,t)H,x(1),(u,v)\}-\{(u,v)H,(z,t),x(1)\}+\{x(1)H,(u,v),(z,t)\}=&\\
&=-\{z(1)-t(2),x(1),(u,v)\}-\{u(1)-v(2),(z,t),x(1)\}+\{x(1),(u,v),(z,t)\}=&\\
&=-z(1)x(1)(u,v)+z(1)(u,v)x(1)-2z(1)\cdot x(1)(u,v)+&\\
&+t(2)x(1)(u,v)-t(2)(u,v)x(1)+2t(2)\cdot x(1)(u,v)-&\\
&-u(1)(z,t)x(1)+u(1)x(1)(z,t)-2u(1)\cdot (z,t)x(1)+&\\
&+v(2)(z,t)x(1)-v(2)x(1)(z,t)+2v(2)\cdot (z,t)x(1)+&\\
&+x(1)(u,v)(z,t)-x(1)(z,t)(u,v)+2x(1)\cdot (u,v)(z,t)=&\\
&=-F\langle z,x\rangle(u,v)+(F\langle z,u\rangle+\frac{1}{2}H\langle z,v\rangle)x(1)-2z(1)(F\langle x,u\rangle+\frac{1}{2}H\langle x,v\rangle)+&\\
&+\frac{1}{2}H\langle t,x\rangle(u,v)-(\frac{1}{2}H\langle t,u\rangle+E\langle t,v\rangle)x(1)+2t(2)(F\langle x,u\rangle+\frac{1}{2}H\langle x,v\rangle)-&\\
&-(F\langle u,z\rangle+\frac{1}{2}H\langle u,t\rangle)x(1)+F\langle u,x\rangle(z,t)-2u(1)(F\langle z,x\rangle+\frac{1}{2}H\langle t,x\rangle)+&\\
&+(\frac{1}{2}H\langle v,z\rangle+E\langle v,t\rangle)x(1)-\frac{1}{2}H\langle v,x\rangle(z,t)+2v(2)(F\langle z,x\rangle+\frac{1}{2}H\langle t,x\rangle)+&\\
&+(F\langle x,u\rangle+\frac{1}{2}H\langle x,v\rangle)(z,t)-(F\langle x,z\rangle+\frac{1}{2}H\langle x,t\rangle)(u,v)+&\\
&+2x(1)(F\langle u,z\rangle+\frac{1}{2}H(\langle u,t\rangle+\langle v,z\rangle)+E\langle v,t\rangle)=&
\ees
\bes
&=-v(1)\langle z,x\rangle-\frac{1}{2}x(1)\langle z,v\rangle-z(1)\langle x,v\rangle+&\\
&+\frac{1}{2}(-u(1)+v(2))\langle t,x\rangle+\frac{1}{2}x(1)\langle t,u\rangle+x(2)\langle t,v\rangle-2t(1)\langle x,u\rangle-t(2)\langle x,v\rangle+&\\
&+\frac{1}{2}x(1)\langle u,t\rangle+t(1)\langle u,x\rangle-u(1)\langle t,x\rangle-&\\
&-\frac{1}{2}x(1)\langle v,z\rangle-x(2)\langle v,t\rangle-\frac{1}{2}(-z(1)+t(2))\langle v,x\rangle-2v(1)\langle z,x\rangle-v(2)\langle t,x\rangle+&\\
&+t(1)\langle x,u\rangle +\frac{1}{2}(-z(1)+t(2))\langle x,v\rangle-v(1)\langle x,z\rangle-\frac{1}{2}(-u(1)+v(2))\langle x,t\rangle+&\\
&+x(1)(\langle u,t\rangle+\langle v,z\rangle)+2x(2)\langle v,t\rangle=&\\
&=-2(u(1)\langle t,x\rangle+t(1)\langle x,u\rangle)-2(v(1)\langle z,x\rangle+z(1)\langle x,v\rangle) + x(1)(\langle u,t\rangle+\langle v,z\rangle)=&\\
&\overset{(\ref{ee28})}{=}3x(1)(\langle u,t\rangle+\langle v,z\rangle)&
\ees
and
\bes
&y(2)\alpha((u,v),(z,t),H)=p(y(2),(u,v),(z,t),H)=&\\
&=-\{(z,t)H,y(2),(u,v)\}-\{(u,v)H,(z,t),y(2)\}+\{y(2)H,(u,v),(z,t)\}=&\\
&=-\{z(1)-t(2),y(2),(u,v)\}-\{u(1)-v(2),(z,t),y(2)\}-\{y(2),(u,v),(z,t)\}=&\\
&=-z(1)y(2)(u,v)+z(1)(u,v)y(2)-2z(1)\cdot y(2)(u,v)+&\\
&+t(2)y(2)(u,v)-t(2)(u,v)y(2)+2t(2)\cdot y(2)(u,v)-&\\
&-u(1)(z,t)y(2)+u(1)y(2)(z,t)-2u(1)\cdot (z,t)y(2)+&\\
&+v(2)(z,t)y(2)-v(2)y(2)(z,t)+2v(2)\cdot (z,t)y(2)-&\\
&-y(2)(u,v)(z,t)+y(2)(z,t)(u,v)-2y(2)\cdot (u,v)(z,t)=&\\
&=-\frac{1}{2}H\langle z,y\rangle(u,v)+(F\langle z,u\rangle+\frac{1}{2}H\langle z,v\rangle)y(2)-2z(1)(\frac{1}{2}H\langle y,u\rangle+E\langle y,v\rangle)+&\\
&+E\langle t,y\rangle(u,v)-(\frac{1}{2}H\langle t,u\rangle+E\langle t,v\rangle)y(2)+2t(2)(\frac{1}{2}H\langle y,u\rangle+E\langle y,v\rangle)-&\\
&-(F\langle u,z\rangle+\frac{1}{2}H\langle u,t\rangle)y(2)+\frac{1}{2}H\langle u,y\rangle(z,t)-2u(1)(\frac{1}{2}H\langle z,y\rangle+E\langle t,y\rangle)+&\\
&+(\frac{1}{2}H\langle v,z\rangle+E\langle v,t\rangle)y(2)-E\langle v,y\rangle(z,t)+2v(2)(\frac{1}{2}H\langle z,y\rangle+E\langle t,y\rangle)-&\\
&-(\frac{1}{2}H\langle y,u\rangle+E\langle y,v\rangle)(z,t)+(\frac{1}{2}H\langle y,z\rangle+E\langle y,t\rangle)(u,v)-&\\
&-2y(2)(F\langle u,z\rangle+\frac{1}{2}H(\langle u,t\rangle+\langle v,z\rangle)+E\langle v,t\rangle)=&\\
&=-\frac{1}{2}(-u(1)+v(2))\langle z,y\rangle+y(1)\langle z,u\rangle+\frac{1}{2}y(2)\langle z,v\rangle-z(1)\langle y,u\rangle-2z(2)\langle y,v\rangle-&\\
&-u(2)\langle t,y\rangle-\frac{1}{2}y(2)\langle t,u\rangle-t(2)\langle y,u\rangle-&\\
&-y(1)\langle u,z\rangle-\frac{1}{2}y(2)\langle u,t\rangle+\frac{1}{2}(-z(1)+t(2))\langle u,y\rangle-u(1)\langle z,y\rangle-2u(2)\langle t,y\rangle+&\\
&+\frac{1}{2}y(2)\langle v,z\rangle+z(2)\langle v,y\rangle-v(2)\langle z,y\rangle-&\\
&-\frac{1}{2}(-z(1)+t(2))\langle y,u\rangle+z(2)\langle y,v\rangle+\frac{1}{2}(-u(1)+v(2))\langle y,z\rangle-u(2)\langle y,t\rangle+&\\
&+2y(1)\langle u,z\rangle+y(2)(\langle u,t\rangle+\langle v,z\rangle)=&\\
&=2(v(2)\langle y,z\rangle+z(2)\langle v,y\rangle)+2(u(2)\langle y,t\rangle+t(2)\langle u,y\rangle)+y(2)(\langle u,t\rangle+\langle v,z\rangle)=&\\
&\overset{(\ref{ee28})}{=}3y(2)(\langle u,t\rangle+\langle v,z\rangle).
\ees
Then it easily follows that $\alpha((u,v),(z,t),E)$, $\alpha((u,v),(z,t),F)$ and $\alpha((u,v),(z,t),H)$ are elements of $\Gamma(\mathcal{M}).$

\begin{center}
    \textbf{The case $\alpha((u,v),(z,t),(p,q))$}
\end{center}
Finally, we have
\bes
&E\alpha((u,v),(z,t),(p,q))=p(E,(u,v),(z,t),(p,q))=&\\
&=-\{(z,t)(p,q),E,(u,v)\}-\{(u,v)(p,q),(z,t),E\}+\{E(p,q),(u,v),(z,t)\}=&\\
&=-\{F\langle z,p\rangle + \frac{1}{2}H(\langle z,q\rangle+\langle t,p\rangle)+E\langle t,q\rangle,E,(u,v)\}-&\\
&-\{F\langle u,p\rangle + \frac{1}{2}H(\langle u,q\rangle+\langle v,p\rangle)+E\langle v,q\rangle,(z,t),E\}-\{p(2),(u,v),(z,t)\}=&\\
&=-((F\langle z,p\rangle + \frac{1}{2}H(\langle z,q\rangle+\langle t,p\rangle)+E\langle t,q\rangle)E (u,v)-&\\
&-(F\langle z,p\rangle + \frac{1}{2}H(\langle z,q\rangle+\langle t,p\rangle)+E\langle t,q\rangle)(u,v)E +&\\
&+2(F\langle z,p\rangle + \frac{1}{2}H(\langle z,q\rangle+\langle t,p\rangle)+E\langle t,q\rangle)\cdot E(u,v))-&\\
&-((F\langle u,p\rangle + \frac{1}{2}H(\langle u,q\rangle+\langle v,p\rangle)+E\langle v,q\rangle)(z,t)E-&\\
&-(F\langle u,p\rangle + \frac{1}{2}H(\langle u,q\rangle+\langle v,p\rangle)+E\langle v,q\rangle)E(z,t)+&\\
&+2(F\langle u,p\rangle + \frac{1}{2}H(\langle u,q\rangle+\langle v,p\rangle)+E\langle v,q\rangle)\cdot (z,t)E)-&\\
&-(p(2)(u,v)(z,t)-p(2)(z,t)(u,v)+2p(2)\cdot (u,v)(z,t))=&\\
&=(\frac{1}{2}H\langle z,p\rangle+\frac{1}{2}E(\langle z,q\rangle + \langle t,p\rangle))(u,v)+&\\
&+(v(1)\langle z,p\rangle-\frac{1}{2}u(1)(\langle z,q\rangle + \langle t,p\rangle)-u(2)\langle t,q\rangle)E+&\\
&+2(F\langle z,p\rangle + \frac{1}{2}H(\langle z,q\rangle+\langle t,p\rangle)+E\langle t,q\rangle)u(2)+&\\
&(-t(1)\langle u,p\rangle+\frac{1}{2}z(1)(\langle u,q\rangle+\langle v,p\rangle)-\frac{1}{2}t(2)(\langle u,q\rangle+\langle v,p\rangle)+z(2)\langle v,q\rangle)E+&\\
&(-\frac{1}{2}H\langle u,p\rangle-\frac{1}{2}E(\langle u,q\rangle+\langle v,p\rangle)(z,t)-&\\
&-2(F\langle u,p\rangle + \frac{1}{2}H(\langle u,q\rangle+\langle v,p\rangle)+E\langle v,q\rangle)z(2)-&\\
&-(\frac{1}{2}H\langle p,u\rangle+E\langle p,v\rangle)(z,t)+(\frac{1}{2}H\langle p,z\rangle+E\langle p,t\rangle)(u,v)-&\\
&-2p(2)(F\langle u,z\rangle+\frac{1}{2}H\langle u,t\rangle+\frac{1}{2}H\langle v,z\rangle+E\langle v,t\rangle)=&\\
&=-\frac{1}{2}u(1)\langle z,p\rangle+\frac{1}{2}v(2)\langle z,p\rangle-\frac{1}{2}u(2)(\langle z,q\rangle+\langle t,q\rangle)+&\\
&+v(2)\langle z,p\rangle-\frac{1}{2}u(2)(\langle z,q\rangle+\langle t,p\rangle)+2u(1)\langle z,p\rangle+u(2)(\langle z,q\rangle+\langle t,p\rangle)-&\\
&-t(2)\langle u,p\rangle+\frac{1}{2}z(2)(\langle u,q\rangle+\langle v,p\rangle)+\frac{1}{2}z(1)\langle u,p\rangle-\frac{1}{2}t(2)\langle u,p\rangle+\frac{1}{2}z(2)(\langle u,q\rangle+\langle v,p\rangle)-&\\
&-2z(1)\langle u,p\rangle-z(2)(\langle u,q\rangle + \langle v,p\rangle)+&\\
&+\frac{1}{2}z(1)\langle p,u\rangle-\frac{1}{2}t(2)\langle p,u\rangle+z(2)\langle p,v\rangle-\frac{1}{2}u(1)\langle p,z\rangle+\frac{1}{2}v(2)\langle p,z\rangle-u(2)\langle p,t\rangle+&\\
&+2p(1)\langle u,z\rangle +p(2)\langle u,t\rangle+p(2)\langle v,z\rangle=&\\
&=2(u(1)\langle z,p\rangle+z(1)\langle p,u\rangle+p(1)\langle u,z\rangle)+\\
&+t(2)\langle p,u\rangle+p(2)\langle u,t\rangle+u(2)\langle t,p\rangle+&\\
&+v(2)\langle z,p\rangle+z(2)\langle p,v\rangle+p(2)\langle v,z\rangle\overset{(\ref{ee28})}{=}0.&
\ees
Also,
\bes
&F\alpha((u,v),(z,t),(p,q))=p(F,(u,v),(z,t),(p,q))=&\\
&=-\{(z,t)(p,q),F,(u,v)\}-\{(u,v)(p,q),(z,t),F\}+\{F(p,q),(u,v),(z,t)\}=&\\
&=-\{F\langle z,p\rangle + \frac{1}{2}H(\langle z,q\rangle+\langle t,p\rangle)+E\langle t,q\rangle,F,(u,v)\}-&\\
&-\{F\langle u,p\rangle + \frac{1}{2}H(\langle u,q\rangle+\langle v,p\rangle)+E\langle v,q\rangle,(z,t),F\}+\{q(1),(u,v),(z,t)\}=&\\
&=-((F\langle z,p\rangle + \frac{1}{2}H(\langle z,q\rangle+\langle t,p\rangle)+E\langle t,q\rangle)F(u,v)-&\\
&-(F\langle z,p\rangle + \frac{1}{2}H(\langle z,q\rangle+\langle t,p\rangle)+E\langle t,q\rangle)(u,v)F+&\\
&+2(F\langle z,p\rangle + \frac{1}{2}H(\langle z,q\rangle+\langle t,p\rangle)+E\langle t,q\rangle)\cdot F(u,v))-&\\
&-((F\langle u,p\rangle + \frac{1}{2}H(\langle u,q\rangle+\langle v,p\rangle)+E\langle v,q\rangle)(z,t)F-&\\
&-(F\langle u,p\rangle + \frac{1}{2}H(\langle u,q\rangle+\langle v,p\rangle)+E\langle v,q\rangle)F(z,t)+&\\
&+2(F\langle u,p\rangle + \frac{1}{2}H(\langle u,q\rangle+\langle v,p\rangle)+E\langle v,q\rangle)\cdot (z,t)F)+&\\
&+q(1)(u,v)(z,t)-q(1)(z,t)(u,v)+2q(1)\cdot (u,v)(z,t)=&\\
&=-(\frac{1}{2}F(\langle z,q\rangle+\langle t,p\rangle)+\frac{1}{2}H\langle t,q\rangle)(u,v)+&\\
&+(v(1)\langle z,p\rangle-\frac{1}{2}u(1)(\langle z,q\rangle+\langle t,p\rangle)+\frac{1}{2}v(2)(\langle z,q\rangle+\langle t,p\rangle)-u(2)\langle t,q\rangle)F-&\\
&-2(F\langle z,p\rangle + \frac{1}{2}H(\langle z,q\rangle+\langle t,p\rangle)+E\langle t,q\rangle)v(1)-&\\
&-(t(1)\langle u,p\rangle-\frac{1}{2}z(1)(\langle u,q\rangle+\langle v,p\rangle)+\frac{1}{2}t(2)(\langle u,q\rangle+\langle v,p\rangle)-z(2)\langle v,q\rangle)F+&\\
&+(\frac{1}{2}F(\langle u,q\rangle+\langle v,p\rangle)+\frac{1}{2}H\langle v,q\rangle)(z,t)+&\\
&+2(F\langle u,p\rangle + \frac{1}{2}H(\langle u,q\rangle+\langle v,p\rangle)+E\langle v,q\rangle)t(1)+&\\
&+(F\langle q,u\rangle+\frac{1}{2}H\langle q,v\rangle)(z,t)-(F\langle q,z\rangle+\frac{1}{2}H\langle q,t\rangle)(u,v)+&\\
&+2q(1)(F\langle u,z\rangle+\frac{1}{2}H\langle u,t\rangle+\frac{1}{2}H\langle v,z\rangle+E\langle v,t\rangle)=&\\
&=-\frac{1}{2}v(1)(\langle z,q\rangle+\langle t,p\rangle)+\frac{1}{2}u(1)\langle t,q\rangle-\frac{1}{2}v(2)\langle t,q\rangle-&\\
&-\frac{1}{2}v(1)(\langle z,q\rangle+\langle t,p\rangle)+u(1)\langle t,q\rangle+&\\
&+v(1)(\langle z,q\rangle+\langle t,p\rangle)+2v(2)\langle t,q\rangle)+&\\
&+\frac{1}{2}t(1)(\langle u,q\rangle+\langle v,p\rangle)-z(1)\langle v,q\rangle+&\\
&+\frac{1}{2}t(1)(\langle u,q\rangle+\langle v,p\rangle)-\frac{1}{2}z(1)\langle v,q\rangle+\frac{1}{2}t(2)\langle v,q\rangle-&\\
&-t(1)(\langle u,q\rangle+\langle v,p\rangle)-2t(2)\langle v,q\rangle+&\\
&+t(1)\langle q,u\rangle-\frac{1}{2}z(1)\langle q,v\rangle+\frac{1}{2}t(2)\langle q,v\rangle-&\\
&-v(1)\langle q,z\rangle+\frac{1}{2}u(1)\langle q,t\rangle-\frac{1}{2}v(2)\langle q,t\rangle+&\\
&+q(1)(\langle u,t\rangle+\langle v,z\rangle)+2q(2)\langle v,t\rangle)=&\\
&=u(1)\langle t,q\rangle+t(1)\langle q,u\rangle+q(1)\langle u,t\rangle+&\\
&+z(1)\langle q,v\rangle+q(1)\langle v,z\rangle+v(1)\langle z,q\rangle+&\\
&+2(v(2)\langle t,q\rangle+t(2)\langle q,v\rangle+q(2)\langle v,t\rangle)\overset{(\ref{ee28})}{=}0.&
\ees
In addition,
\bes
&H\alpha((u,v),(z,t),(p,q))=p(H,(u,v),(z,t),(p,q))=&\\
&=-\{(z,t)(p,q),H,(u,v)\}-\{(u,v)(p,q),(z,t),H\}+\{H(p,q),(u,v),(z,t)\}=&
\\
&=-\{F\langle z,p\rangle + \frac{1}{2}H(\langle z,q\rangle+\langle t,p\rangle)+E\langle t,q\rangle,H,(u,v)\}-&\\
&-\{F\langle u,p\rangle + \frac{1}{2}H(\langle u,q\rangle+\langle v,p\rangle)+E\langle v,q\rangle,(z,t),H\}+\{-p(1) + q(2),(u,v),(z,t)\}=&\\
&=-((F\langle z,p\rangle + \frac{1}{2}H(\langle z,q\rangle+\langle t,p\rangle)+E\langle t,q\rangle)H(u,v)-&\\
&-(F\langle z,p\rangle + \frac{1}{2}H(\langle z,q\rangle+\langle t,p\rangle)+E\langle t,q\rangle)(u,v)H+&\\
&+2(F\langle z,p\rangle + \frac{1}{2}H(\langle z,q\rangle+\langle t,p\rangle)+E\langle t,q\rangle)\cdot H(u,v))-&\\
&-((F\langle u,p\rangle + \frac{1}{2}H(\langle u,q\rangle+\langle v,p\rangle)+E\langle v,q\rangle)(z,t)H-&\\
&-(F\langle u,p\rangle + \frac{1}{2}H(\langle u,q\rangle+\langle v,p\rangle)+E\langle v,q\rangle)H(z,t)+&\\
&+2(F\langle u,p\rangle + \frac{1}{2}H(\langle u,q\rangle+\langle v,p\rangle)+E\langle v,q\rangle)\cdot (z,t)H)-&\\
&-(p(1)(u,v)(z,t)-p(1)(z,t)(u,v)+2p(1)\cdot (u,v)(z,t))+&\\
&+q(2)(u,v)(z,t)-q(2)(z,t)(u,v)+2q(2)\cdot (u,v)(z,t)=&\\
&=-(-F\langle z,p\rangle+E\langle t,q\rangle)(u,v)+&\\
&+(v(1)\langle z,p\rangle-\frac{1}{2}u(1)(\langle z,q\rangle+\langle t,p\rangle)+\frac{1}{2}v(2)(\langle z,q\rangle+\langle t,p\rangle)-u(2)\langle t,q\rangle)H-&\\
&-2(F\langle z,p\rangle + \frac{1}{2}H(\langle z,q\rangle+\langle t,p\rangle)+E\langle t,q\rangle)(-u(1)+v(2))-&\\
&-(t(1)\langle u,p\rangle-\frac{1}{2}z(1)(\langle u,q\rangle+\langle v,p\rangle)+\frac{1}{2}t(2)(\langle u,q\rangle+\langle v,p\rangle)-z(2)\langle v,q\rangle)H+&\\
&+(-F\langle u,p\rangle+E\langle v,q\rangle)(z,t)-&\\
&-2(F\langle u,p\rangle + \frac{1}{2}H(\langle u,q\rangle+\langle v,p\rangle)+E\langle v,q\rangle)(z(1)-t(2))-&\\
&-(F\langle p,u\rangle+\frac{1}{2}H\langle p,v\rangle)(z,t)+(F\langle p,z\rangle+\frac{1}{2}H\langle p,t\rangle)(u,v)-&\\
&-2p(1)(F\langle u,z\rangle+\frac{1}{2}H\langle u,t\rangle+\frac{1}{2}H\langle v,z\rangle+E\langle v,t\rangle)+&\\
&+(\frac{1}{2}H\langle q,u\rangle+E\langle q,v\rangle)(z,t)-(\frac{1}{2}H\langle q,z\rangle+E\langle q,t\rangle)(u,v)+&\\
&+2q(2)(F\langle u,z\rangle+\frac{1}{2}H\langle u,t\rangle+\frac{1}{2}H\langle v,z\rangle+E\langle v,t\rangle)=&\\
&=v(1)\langle z,p\rangle+u(2)\langle t,q\rangle+&\\
&+v(1)\langle z,p\rangle-\frac{1}{2}u(1)(\langle z,q\rangle+\langle t,p\rangle)-\frac{1}{2}v(2)(\langle z,q\rangle+\langle t,p\rangle)+u(2)\langle t,q\rangle-&\\
&-2v(1)\langle z,p\rangle-u(1)(\langle z,q\rangle+\langle t,p\rangle)-v(2)(\langle z,q\rangle+\langle t,p\rangle)-2u(2)\langle t,q\rangle-&\\
&-t(1)\langle u,p\rangle+\frac{1}{2}z(1)(\langle u,q\rangle+\langle v,p\rangle)+\frac{1}{2}t(2)(\langle u,q\rangle+\langle v,p\rangle)-z(2)\langle v,q\rangle-&\\
&-t(1)\langle u,p\rangle-z(2)\langle v,q\rangle -&\\
&2t(1)\langle u,p\rangle+z(1)(\langle u,q\rangle+\langle v,p\rangle)+t(2)(\langle u,q\rangle+\langle v,p\rangle)+2z(2)\langle v,q\rangle-&\\
&-t(1)\langle p,u\rangle+\frac{1}{2}z(1)\langle p,v\rangle-\frac{1}{2}t(2)\langle p,v\rangle+v(1)\langle p,z\rangle-\frac{1}{2}u(1)\langle p,t\rangle+\frac{1}{2}v(2)\langle p,t\rangle-&\\
&-p(1)\langle u,t\rangle-p(1)\langle v,z\rangle-2p(2)\langle v,t\rangle+&\\
&-\frac{1}{2}z(1)\langle q,u\rangle+\frac{1}{2}t(2)\langle q,u\rangle-z(2)\langle q,v\rangle+\frac{1}{2}u(1)\langle q,z\rangle-\frac{1}{2}v(2)\langle q,z\rangle+u(2)\langle q,t\rangle+&\\
&-2q(1)\langle u,z\rangle-q(2)\langle u,t\rangle-q(2)\langle v,z\rangle=&
\ees
\bes
&=2(u(1)\langle q,z\rangle+q(1)\langle z,u\rangle+z(1)\langle u,q\rangle)+&\\
&+2(v(2)\langle p,t\rangle+p(2)\langle t,v\rangle+t(2)\langle v,p\rangle)+&\\
&+u(1)\langle p,t\rangle + p(1)\langle t,u\rangle+ t(1)\langle u,p\rangle +&\\
&+z(1)\langle v,p\rangle + v(1)\langle p,z\rangle + p(1)\langle z,v\rangle+&\\
&+z(2)\langle v,q\rangle + v(2)\langle q,z\rangle + q(2)\langle z,v\rangle+&\\
&+u(2)\langle q,t\rangle+q(2)\langle t,u\rangle+t(2)\langle u,q\rangle\overset{(\ref{ee28})}{=}0.&
\ees
Besides, 
\bes
&x(1)\alpha((u,v),(z,t),(p,q))=p(x(1),(u,v),(z,t),(p,q))=&\\
&=-\{(z,t)(p,q),x(1),(u,v)\}-\{(u,v)(p,q),(z,t),x(1)\}+\{x(1)(p,q),(u,v),(z,t)\}=&\\
&=-\{F\langle z,p\rangle + \frac{1}{2}H(\langle z,q\rangle+\langle t,p\rangle)+E\langle t,q\rangle,x(1),(u,v)\}-&\\
&-\{F\langle u,p\rangle + \frac{1}{2}H(\langle u,q\rangle+\langle v,p\rangle)+E\langle v,q\rangle,(z,t),x(1)\}+&\\
&+\{F\langle x,p\rangle+\frac{1}{2}H\langle x,q\rangle,(u,v),(z,t)\}=&\\
&=-((F\langle z,p\rangle + \frac{1}{2}H(\langle z,q\rangle+\langle t,p\rangle)+E\langle t,q\rangle)x(1)(u,v)-&\\
&-(F\langle z,p\rangle + \frac{1}{2}H(\langle z,q\rangle+\langle t,p\rangle)+E\langle t,q\rangle)(u,v)x(1)+&\\
&+2(F\langle z,p\rangle + \frac{1}{2}H(\langle z,q\rangle+\langle t,p\rangle)+E\langle t,q\rangle)\cdot x(1)(u,v))-&\\
&-((F\langle u,p\rangle + \frac{1}{2}H(\langle u,q\rangle+\langle v,p\rangle)+E\langle v,q\rangle)(z,t)x(1)-&\\
&-(F\langle u,p\rangle + \frac{1}{2}H(\langle u,q\rangle+\langle v,p\rangle)+E\langle v,q\rangle)x(1)(z,t)+&\\
&+2(F\langle u,p\rangle + \frac{1}{2}H(\langle u,q\rangle+\langle v,p\rangle)+E\langle v,q\rangle)\cdot (z,t)x(1))+&\\
&+\{F,(u,v),(z,t)\}\langle x,p\rangle+\frac{1}{2}\{H,(u,v),(z,t)\}\langle x,q\rangle=&\\
&\overset{(\ref{lie1})-(\ref{lie4})}{=}(\frac{1}{2}x(1)(\langle z,q\rangle+\langle t,p\rangle)+x(2)\langle t,q\rangle)(u,v)+&\\
&+(v(1)\langle z,p\rangle-\frac{1}{2}(u(1)-v(2))(\langle z,q\rangle+\langle t,p\rangle)-u(2)\langle t,q\rangle)x(1)-&\\
&-2(F\langle z,p\rangle + \frac{1}{2}H(\langle z,q\rangle+\langle t,p\rangle)+E\langle t,q\rangle)(F\langle x,u\rangle+\frac{1}{2}H\langle x,v\rangle)-&\\
&-(t(1)\langle u,p\rangle-\frac{1}{2}(z(1)-t(2))(\langle u,q\rangle+\langle v,p\rangle)-z(2)\langle v,q\rangle)x(1)-&\\
&-(\frac{1}{2}x(1)(\langle u,q\rangle+\langle v,p\rangle)+x(2)\langle v,q\rangle)(z,t)-&\\
&-2(F\langle u,p\rangle + \frac{1}{2}H(\langle u,q\rangle+\langle v,p\rangle)+E\langle v,q\rangle)(F\langle z,x\rangle+\frac{1}{2}H\langle t,x\rangle)=&\\
&=\frac{1}{2}F\langle x,u\rangle(\langle z,q\rangle+\langle t,p\rangle)+\frac{1}{4}H\langle x,v\rangle(\langle z,q\rangle+\langle t,p\rangle)+(\frac{1}{2}H\langle x,u\rangle+E\langle x,v\rangle)\langle t,q\rangle+&\\
&+F\langle v,x\rangle\langle z,p\rangle-\frac{1}{2}(F\langle u,x\rangle-\frac{1}{2}H\langle v,x\rangle)(\langle z,q\rangle+\langle t,p\rangle)-\frac{1}{2}H\langle u,x\rangle\langle t,q\rangle+&\\
&+F\langle z,p\rangle\langle x,v\rangle-F\langle x,u\rangle(\langle z,q\rangle+\langle t,p\rangle)-H\langle t,q\rangle\langle x,u\rangle-E\langle t,q\rangle\langle x,v\rangle-&\\
&-F\langle t,x\rangle\langle u,p\rangle+\frac{1}{2}(F\langle z,x\rangle-\frac{1}{2}H\langle t,x\rangle)(\langle u,q\rangle+\langle v,p\rangle)+\frac{1}{2}H\langle z,x\rangle\langle v,q\rangle-&\\
&-\frac{1}{2}(F\langle x,z\rangle+\frac{1}{2}H\langle x,t\rangle)(\langle u,q\rangle+\langle v,p\rangle)-\frac{1}{2}H\langle x,z\rangle\langle v,q\rangle-E\langle x,t\rangle\langle v,q\rangle+&\\
&+F\langle u,p\rangle\langle t,x\rangle-F(\langle u,q\rangle+\langle v,p\rangle)\langle z,x\rangle-H\langle v,q\rangle\langle z,x\rangle-E\langle v,q\rangle\langle t,x\rangle=0&
\ees
and
\bes
&y(2)\alpha((u,v),(z,t),(p,q))=p(y(2),(u,v),(z,t),(p,q))=&\\
&=-\{(z,t)(p,q),y(2),(u,v)\}-\{(u,v)(p,q),(z,t),y(2)\}+\{y(2)(p,q),(u,v),(z,t)\}=&\\
&=-\{F\langle z,p\rangle + \frac{1}{2}H(\langle z,q\rangle+\langle t,p\rangle)+E\langle t,q\rangle,y(2),(u,v)\}-&\\
&-\{F\langle u,p\rangle + \frac{1}{2}H(\langle u,q\rangle+\langle v,p\rangle)+E\langle v,q\rangle,(z,t),y(2)\}+&\\
&+\{\frac{1}{2}H\langle y,p\rangle+E\langle y,q\rangle,(u,v),(z,t)\}=&\\
&=-((F\langle z,p\rangle + \frac{1}{2}H(\langle z,q\rangle+\langle t,p\rangle)+E\langle t,q\rangle)y(2)(u,v)-&\\
&-(F\langle z,p\rangle + \frac{1}{2}H(\langle z,q\rangle+\langle t,p\rangle)+E\langle t,q\rangle)(u,v)y(2)+&\\
&+2(F\langle z,p\rangle + \frac{1}{2}H(\langle z,q\rangle+\langle t,p\rangle)+E\langle t,q\rangle)\cdot y(2)(u,v))-&\\
&-((F\langle u,p\rangle + \frac{1}{2}H(\langle u,q\rangle+\langle v,p\rangle)+E\langle v,q\rangle)(z,t)y(2)-&\\
&-(F\langle u,p\rangle + \frac{1}{2}H(\langle u,q\rangle+\langle v,p\rangle)+E\langle v,q\rangle)y(2)(z,t)+&\\
&+2(F\langle u,p\rangle + \frac{1}{2}H(\langle u,q\rangle+\langle v,p\rangle)+E\langle v,q\rangle)\cdot (z,t)y(2))+&\\
&+\frac{1}{2}\{H,(u,v),(z,t)\}\langle y,p\rangle+\{E,(u,v),(z,t)\}\langle y,q\rangle=&\\
&\overset{(\ref{lie1})-(\ref{lie4})}{=}-(y(1)\langle z,p\rangle+\frac{1}{2}y(2)(\langle z,q\rangle+\langle t,p\rangle))(u,v)+&\\
&+(v(1)\langle z,p\rangle-\frac{1}{2}(u(1)-v(2))(\langle z,q\rangle+\langle t,p\rangle)-u(2)\langle t,q\rangle)y(2)-&\\
&-2(F\langle z,p\rangle + \frac{1}{2}H(\langle z,q\rangle+\langle t,p\rangle)+E\langle t,q\rangle)(\frac{1}{2}H\langle y,u\rangle+E\langle y,v\rangle)-&\\
&-(t(1)\langle u,p\rangle-\frac{1}{2}(z(1)-t(2))(\langle u,q\rangle+\langle v,p\rangle)-z(2)\langle v,q\rangle)y(2)+&\\
&+(y(1)\langle u,p\rangle+\frac{1}{2}y(2)(\langle u,q\rangle+\langle v,p\rangle))(z,t)-&\\
&-2(F\langle u,p\rangle + \frac{1}{2}H(\langle u,q\rangle+\langle v,p\rangle)+E\langle v,q\rangle)(\frac{1}{2}H\langle z,y\rangle+E\langle t,y\rangle)=&\\
&=-(F\langle y,u\rangle+\frac{1}{2}H\langle y,v\rangle)\langle z,p\rangle-\frac{1}{2}(\frac{1}{2}H\langle y,u\rangle+E\langle y,v\rangle)(\langle z,q\rangle+\langle t,p\rangle))+&\\
&+\frac{1}{2}H\langle v,y\rangle\langle z,p\rangle-\frac{1}{2}(\frac{1}{2}H\langle u,y\rangle-E\langle v,y\rangle)(\langle z,q\rangle+\langle t,p\rangle)-E\langle u,y\rangle\langle t,q\rangle+&\\
&+F\langle z,p\rangle\langle y,u\rangle+H\langle z,p\rangle\langle y,v\rangle+E(\langle z,q\rangle+\langle t,p\rangle)\langle y,v\rangle-E\langle t,q\rangle\langle y,u\rangle-&\\
&-\frac{1}{2}H\langle t,y\rangle\langle u,p\rangle+\frac{1}{2}(\frac{1}{2}H\langle z,y\rangle-E\langle t,y\rangle)(\langle u,q\rangle+\langle v,p\rangle)+E\langle z,y\rangle\langle v,q\rangle+&\\
&+(F\langle y,z\rangle+\frac{1}{2}H\langle y,t\rangle)\langle u,p\rangle+\frac{1}{2}(\frac{1}{2}H\langle y,z\rangle+E\langle y,t\rangle)(\langle u,q\rangle+\langle v,p\rangle)+&\\
&+F\langle u,p\rangle\langle z,y\rangle+H\langle u,p\rangle\langle t,y\rangle+E(\langle u,q\rangle+\langle v,p\rangle)\langle t,y\rangle-E\langle v,q\rangle\langle z,y\rangle=0.&
\ees
Thus $\alpha((u,v),(z,t),(p,q))|_\mathcal{M}=0;$ so $\alpha((u,v),(z,t),(p,q))$ belongs to $\Gamma(\mathcal{M}).$

Consequently, from the above cases we get that $\mathcal{M}$ belongs to the variety $\mathcal{H}.$ Therefore, the theorem is proved. 
\end{dem}

We know that the commutator in any alternative algebra satisfies the identities $(\ref{e})$, and so every Lie algebra is a Malcev algebra. The speciality problem for Malcev algebras asks if any Malcev algebra is isomorphic to a subalgebra
of the commutator algebra of some alternative algebra. In this case, we have:

\begin{coro}
    Every Malcev algebra $\mathcal{M}$ in $\mathcal{H}$ containing $L\cong\mathfrak{sl}_{2}$ with $mL\neq 0$ for any $0\neq m\in\mathcal{M}$ is special. 
\end{coro}

\subsection{Examples}
\subsubsection{Exceptional Malcev algebras}

Let $\mathcal{B}$ be a unital associative commutative algebra over a field of characteristic $\neq 2,3$  and 
consider the full $2\times 2$ matrix algebra matrices $M_2(\mathcal{B})$ over $\mathcal{B}$  with the following basis:
$$
I=\left(\begin{array}{cc}
\frac{1}{2} & 0 \\
0 & \frac{1}{2}
\end{array}\right),\,\,\, 
E=\left(\begin{array}{cc}
0 & -\frac{1}{2} \\
0 & \,\,\,\,0
\end{array}\right), \,\,\,
H=\left(\begin{array}{cc}
-\frac{1}{2} & 0 \\
\,\,\,\,0 & \frac{1}{2}
\end{array}\right),\,\,\,
F=\left(\begin{array}{cc}
0 & 0 \\
\frac{1}{2} & 0
\end{array}\right).$$

Let $\mathcal{M}=\mathcal{M}_{7}(\mathcal{B})$ be the exceptional Malcev algebra over $\mathcal{B}$. In this case, 
\begin{center}
    $\mathcal{M}=\mathfrak{sl}_{2}(\mathcal{B}) \oplus v M_2(\mathcal{B})$
\end{center}
where $N_{\mathcal{M}}=\mathfrak{sl}_{2}(\mathcal{B})$, $J_{\mathcal{M}}=v M_2(\mathcal{B})$, 
with the product defined by
\begin{equation}\label{relations}
X \cdot Y=XY-YX, \quad X \cdot v A=v\left(X^* A-XA\right), \quad v A \cdot v B=B A^* - A B^*,
\end{equation}
where $X,Y \in \mathfrak{sl}_{2}(\mathcal{B})$, $A, B \in M_2(\mathcal{B})$ and 
$A \mapsto A^*$ is the symplectic involution:
\begin{center}
$\left(\begin{array}{cc}
a & b \\
c & d
\end{array}\right)^*
\mapsto
\left(\begin{array}{cc}
d & -b \\
-c & a
\end{array}\right)$.
\end{center}

Take $V=\mathcal{B}^2=\{(a, b) \mid a, b \in \mathcal{B}\},\,(a, b)(1)=v\left(\begin{array}{cc}0 & 0 \\ a & b\end{array}\right),\,(a, b)(2)=v\left(\begin{array}{cc}-a & -b \\ \,\,\,\,0 & \,\,\,\,0\end{array}\right)$. 
Then using $(\ref{relations})$ it is easy to see that $\{(a,b)(1), (a,b)(2)\}$ satisfies the relations $(\ref{ee35}).$
Hence
we have $\mathcal{M}=\mathfrak{sl}_{2}(\mathcal{B})\oplus V(1)\oplus V(2),$ with $\langle(a, b),(c, d)\rangle=-4\textup{det}\left(\begin{array}{cc}
a & b \\
c & d
\end{array}\right).$

In fact, by $(\ref{ee32})$
$$
\begin{aligned}
F\langle(a, b),(c, d)\rangle & =(a, b)(1) \cdot (c, d)(1) \\
& =v\left(\begin{array}{cc}
0 & 0 \\
a & b
\end{array}\right) \cdot v\left(\begin{array}{cc}
0 & 0 \\
c & d
\end{array}\right) \\
& \overset{(\ref{relations})}{=}\left(\begin{array}{cc}
0 & 0 \\
c & d
\end{array}\right) \left(\begin{array}{cc}
b & 0 \\
-a & 0
\end{array}\right)-\left(\begin{array}{cc}
0 & 0 \\
a & b
\end{array}\right) \left(\begin{array}{ll}
d & 0 \\
-c & 0
\end{array}\right)=-4F\operatorname{det}\left(\begin{array}{ll}
a & b \\
c & d
\end{array}\right).
\end{aligned}
$$
So 
$F(\langle(a, b),(c, d)\rangle +4\operatorname{det}\left(\begin{array}{ll}
a & b \\
c & d
\end{array}\right))=0$ implies
$L(\langle(a, b),(c, d)\rangle +4\operatorname{det}\left(\begin{array}{ll}
a & b \\
c & d
\end{array}\right))=0.$ Thus, 
\begin{center}
    $\langle(a, b),(c, d)\rangle=-4\textup{det}\left(\begin{array}{cc}
a & b \\
c & d
\end{array}\right).$
\end{center}
Now for any $u=(a, b), v=(c, d), w=(e, f) \in V$ we have
$$
\begin{aligned}
w\langle u,v\rangle & + u\langle v,w\rangle + v\langle w,u\rangle  \\
& =-4(e, f)\operatorname{det}\left(\begin{array}{ll}
a & b \\
c & d
\end{array}\right)-4(a, b)\operatorname{det}\left(\begin{array}{cc}
c & d \\
e & f
\end{array}\right)-4(c, d)\operatorname{det}\left(\begin{array}{cc}
e & f \\
a & b
\end{array}\right) \\
& =-4(e, f)(a d-b c)-4(a, b)(c f-d e)-4(c, d)(e b-f a)=(0,0);
\end{aligned}
$$
hence $\mathcal{M}_{7}(\mathcal{B})$ satisfies $(\ref{ee28})$.

The following proposition gives a specific certain condition to obtain a coordinatization theorem for Malcev algebras containing $L=\mathfrak{sl}_{2}(\mathbb{F}).$ 

\begin{pro}\label{p2}
    The algebra $\mathcal{M}=\mathfrak{sl}_{2}(\mathcal{B}) \oplus V^2$ from Theorem \ref{principal} is isomorphic
to the exceptional Malcev algebra $\mathcal{M}_{7}(\mathcal{B})$ if and only if there exist $u, v\in V$ such that $\langle u, v\rangle=1.$
\end{pro}
\begin{dem}
We have already checked that the algebra $\mathcal{M}_{7}(\mathcal{B})$ has the form $\mathfrak{sl}_{2}(\mathcal{B}) \oplus V^2$. It is observed that $\langle u, v\rangle=1$ for $u=(\frac{1}{2},0),$ $v=(0,-\frac{1}{2}) \in V.$

Let now $\mathcal{A}=\mathfrak{sl}_{2}(\mathbb{F})+V_{2}(u)+V_{2}(v)$ be such that there exist $u, v\in V$ with $\langle u, v\rangle=1.$ It follows from Proposition $\ref{p3}$ and its proof that $\mathcal{A}$ is a subalgebra of $\mathcal{M}$ isomorphic to the 7-dimensional exceptional Malcev algebra $\mathcal{M}_{7}(\mathbb{F})$ with $m\algA\neq 0$ for any $0\neq m\in\mathcal{M}$ because $mL\neq 0$ for any $0\neq m\in\mathcal{M}$. Therefore, by \cite{LSS0} (Theorem 3.1), 
    $\mathcal{M} \cong \mathcal{M}_{7}(\mathcal{U})$
for a certain associative commutative algebra $\mathcal{U}$. It follows from $(\ref{ee28})$ that $V=\mathcal{B} \cdot V_{2}(u)+\mathcal{B} \cdot V_{2}(v)$ and $\mathcal{U}=\mathcal{B}$.
Thus, the proposition is proved.
\end{dem}

\subsubsection{Algebras obtained by (commutative) Cayley–Dickson process}

Note that if the mapping $\langle\cdot , \cdot\rangle: V^2 \longrightarrow \mathcal{B}$ is trivial, then the algebra $\mathcal{M}$ is just a split null extension of the algebra $\mathfrak{sl}_{2}(\mathcal{B})$ by a bimodule $V^2$. In this case, $V$ may be an arbitrary commutative $\mathcal{B}$-bimodule. For instance, when $\mathcal{B}=\mathbb{F}$ and $V=\mathbb{F}$ we get in this way the algebra $\mathcal{M}_{5}(\mathbb{F})=\mathfrak{sl}_{2}(\mathbb{F}) \oplus V_{2}(u)$ which is a non-Lie Malcev algebra of dimension five.


\begin{ex}\label{ej8}
The 5-dimensional non-Lie Malcev algebra $\mathcal{M}_{5}(\mathbb{F})=\mathfrak{sl}_{2}(\mathbb{F})\oplus V_2(u)$ that contains $\mathfrak{sl}_{2}(\mathbb{F})$ has the following multiplication table
\begin{center}
\begin{tabular}{|c||c|c|c|c|c|}
\hline
$.$ & $E$ & $F$ & $H$ & $u(1)$ & $u(2)$ \\\hline\hline
$E$ & $0$ & $H/2$ & $E$ & $-u(2)$ & $0$ \\\hline
$F$ & $-H/2$ & $0$ & $-F$ & $0$ & $u(1)$\\\hline
$H$ & $-E$ & $F$ & $0$ & $-u(1)$ & $u(2)$ \\\hline
$u(1)$ & $u(2)$ & $0$ & $u(1)$ & $0$ & $0$ \\\hline
$u(2)$ & $0$ & $-u(1)$ & $-u(2)$ & $0$ & $0$  \\\hline
\end{tabular}
\end{center}
where $\{u(1), u(2)\}$ is the basis of $V_{2}(u).$ 
\end{ex}

If the mapping $\langle\cdot , \cdot\rangle: V^2 \rightarrow \mathcal{B}$ is not trivial, then by $(\ref{ee28})$ the rank of $V$ as a $\mathcal{B}$-bimodule is less than 3. Observe that the left side of $(\ref{ee28})$ is $\mathcal{B}$-multilinear and skew-symmetric on $u, v, w$. Therefore, it holds when $\Lambda^3\left(V_{\mathcal{B}}\right)=0$. In particular it holds if the rank of $V$ is less or equal to 2 . If $V \subseteq \mathcal{B} \cdot x$ then the mapping $\langle\cdot , \cdot\rangle$ is trivial by skew-symmetry. Let us consider now the case when $V$ is a 2-generated $\mathcal{B}$-module.

Let $\mathcal{U}$ be a unital associative commutative algebra over a field of characteristic $\neq 2,3$ and $\alpha \in \mathcal{U}$. Denote by $\widetilde{\mathcal{M}}\left(\mathfrak{sl}_{2}(\mathcal{U}), \alpha\right)$ the algebra $\mathfrak{sl}_{2}(\mathcal{U}) \oplus v M_2(\mathcal{U})$ with a product defined by the following analogue of $(\ref{relations})$:
\begin{equation}\label{relations1}
X \cdot Y=XY-YX, \quad X \cdot v A=v\left(X^* A-XA\right), \quad v A \cdot v B=\alpha (B A^* - A B^*),
\end{equation}
where $X,Y \in \mathfrak{sl}_{2}(\mathcal{B})$, $A, B \in M_2(\mathcal{B})$ and 
$A \mapsto A^*$ is the symplectic involution.

The algebra $\widetilde{\mathcal{M}}\left(\mathfrak{sl}_{2}(\mathcal{U}), \alpha\right)$ is a Malcev algebra containing $\mathfrak{sl}_{2}(\mathcal{U})$. Here
\begin{center}
    $\widetilde{\mathcal{M}}\left(\mathfrak{sl}_{2}(\mathcal{U})), \alpha\right)
    =\textup{CD}\left(M_2(\mathcal{B}), \alpha\right)^{(-)}/\mathbb{F}\cdot 2I
    $
\end{center}
where $2I$ is a matriz identity of $M_2(\mathcal{U}).$
The algebra $\textup{CD}\left(M_2(\mathcal{U}), \alpha\right)=M_{2}(\mathcal{U}) \oplus v M_2(\mathcal{U})$ is an alternative algebra containing $M_2(\mathcal{U})$ with the same identity element, which by \cite{LSS1} we will call it the \textit{algebra obtained from $M_2(\mathcal{U})$ by the Cayley–Dickson process} with parameter $\alpha$.  
The algebra $\widetilde{\mathcal{M}}\left(\mathfrak{sl}_{2}(\mathcal{U}), \alpha\right)$ is an exceptional Malcev algebra if and only if the parameter $\alpha$ is invertible in $\mathcal{U}$.

\begin{theo}\label{t2}
Let $\mathcal{B}$ be a unital associative commutative algebra, let $V=\mathcal{B}^2$, and let $\langle\cdot , \cdot\rangle: V^2 \rightarrow \mathcal{B}$ be a skew-symmetric $\mathcal{B}$-bilinear mapping. Then the algebra $\mathcal{M}=\mathfrak{sl}_{2}(\mathcal{B}) \oplus V^2$ is isomorphic to an algebra $\widetilde{\mathcal{M}}\left(\mathfrak{sl}_{2}(\mathcal{B}), \alpha\right)$, where $\alpha=-\langle(\frac{1}{2},0),(0,\frac{1}{2})\rangle$. Conversely, every algebra $\widetilde{\mathcal{M}}\left(\mathfrak{sl}_{2}(\mathcal{U}), \alpha\right)$ has this form.
\end{theo}
\begin{dem}
Let $\widetilde{\mathcal{M}}=\widetilde{\mathcal{M}}\left(\mathfrak{sl}_{2}(\mathcal{U}), \alpha\right)$. Take $V=\mathcal{U}^2=\{(a, b) \mid a, b \in \mathcal{U}\},(a, b)(1)=v\left(\begin{array}{cc}0 & 0 \\ 
a & b\end{array}\right)$, and 
$(a, b)(2)=v\left(\begin{array}{cc}-a & -b \\ \,\,\,\,0 & \,\,\,\,0\end{array}\right) \in v M_2(\mathcal{U})$. Then we have, as before, $\widetilde{\mathcal{M}}=\mathfrak{sl}_{2}(\mathcal{U}) \oplus V(1) \oplus V(2)$, with $\langle(a, b),(c, d)\rangle=-4\alpha \operatorname{det}\left(\begin{array}{ll}a & b \\ c & d\end{array}\right)$. In particular, $\langle(\frac{1}{2},0),(0,\frac{1}{2})\rangle=-\alpha$.

Conversely, let $\mathcal{M}=\mathfrak{sl}_{2}(\mathcal{B}) \oplus V^2$, where $V \cong \mathcal{B}^2$ and $\langle(\frac{1}{2},0),(0,\frac{1}{2})\rangle=-\alpha$. Define the mapping $\varphi: V^2=V(1) \oplus V(2) \rightarrow v M_2(\mathcal{B}) \subset \widetilde{\mathcal{M}}\left(\mathfrak{sl}_{2}(\mathcal{B}), \alpha\right)$ by sending, for any $a, b \in \mathcal{B}.$

$$
(a, b)(1) \mapsto v\left(\begin{array}{cc}
0 & 0 \\
a & b
\end{array}\right), \quad(a, b)(2) \mapsto v\left(\begin{array}{cc}
-a & -b \\
\,\,\,\,0 & \,\,\,\,0
\end{array}\right)
$$
It is easy to see that $\varphi$ is an isomorphism of Malcev $\mathfrak{sl}_{2}(\mathcal{B})$-modules. Moreover, let $x=(a, b),$ $ y=(c, d) \in V=\mathcal{B}^2$, then we have
$$
\begin{aligned}
\langle x, y\rangle  =\langle(a, b),(c, d)\rangle &=\langle a(1,0)+b(0,1), c(1,0)+d(0,1)\rangle \\
& =(a d-b c)\langle(1,0),(0,1)\rangle=-4\alpha(a d-b c) .
\end{aligned}
$$
Let $z=(e, f),$ $ t=(g, h) \in V$; then we have by $(\ref{ee33})$,
$$
\begin{aligned}
(x, y)(z, t) & = \dfrac{1}{2}\left(\begin{array}{cc}
-\frac{1}{2}(\langle x, t\rangle+\langle y, z\rangle)  & -\langle y,t\rangle \\
\langle x,z\rangle & \frac{1}{2}(\langle x, t\rangle+\langle y, z\rangle)
\end{array}\right)\\
& =\left(\begin{array}{cc}
\alpha(ah-bg+cf-de) & 2\alpha(ch-dg) \\
-2\alpha(af-be) & -\alpha(ah-bg+cf-de) 
\end{array}\right).
\end{aligned}
$$
On the other hand, 
$$
\begin{aligned}
\varphi(x, y) \cdot \varphi(z, t) & =v\left(\begin{array}{cc}
a & b \\
c & d
\end{array}\right) \cdot v\left(\begin{array}{cc}
e & f \\
g & h
\end{array}\right)\\
&\overset{(\ref{relations1})}{=}\alpha[\left(\begin{array}{cc}
e & f \\
g & h
\end{array}\right) \cdot\left(\begin{array}{ll}
d & -b \\
-c & a
\end{array}\right)-
\left(\begin{array}{cc}
a & b \\
c & d
\end{array}\right) \cdot\left(\begin{array}{ll}
h & -f \\
-g & e
\end{array}\right)] \\
& =\alpha\left(\begin{array}{cc}
-bg+ah+cf-de & 2(ch-dg) \\
2(-af+be) & de-cf-ah+bg
\end{array}\right).
\end{aligned}
$$
Therefore, the mapping
$$
\text {id}+\varphi: \mathcal{M}=\mathfrak{sl}_{2}(\mathcal{B}) \oplus V^2 \longrightarrow \widetilde{\mathcal{M}}\left(\mathfrak{sl}_{2}(\mathcal{U}), \alpha\right)=\mathfrak{sl}_{2}(\mathcal{B}) \oplus v M_2(\mathcal{B})
$$
is an isomorphism. 

The theorem is proved.
\end{dem}

\begin{ex}\label{ej9}
We have the 7-dimensional exceptional Malcev algebra $\mathcal{M}_{7}(\mathbb{F})=\mathfrak{sl}_{2}(\mathbb{F})\oplus V_{2}(u)\oplus V_{2}(v)$ that contains $\mathfrak{sl}_{2}(\mathbb{F}).$ The multiplication table of $\mathcal{M}_{7}(\mathbb{F})$ is 
\begin{center}
\begin{tabular}{|c||c|c|c|c|c|c|c|}
\hline
$.$ & $E$ & $F$ & $H$ & $u(1)$ & $u(2)$ & $v(1)$ & $v(2)$\\\hline\hline
$E$ & $0$ & $H/2$ & $E$ & $-u(2)$ & $0$ & $-v(2)$ & $0$\\\hline
$F$ & $-H/2$ & $0$ & $-F$ & $0$ & $u(1)$ & $0$ & $v(1)$\\\hline
$H$ & $-E$ & $F$ & $0$ & $-u(1)$ & $u(2)$  & $-v(1)$ & $v(2)$\\\hline
$u(1)$ & $u(2)$ & $0$ & $u(1)$ & $0$ & $0$ & $F\langle u,v\rangle$ & $H\langle u,v\rangle/2$\\\hline 
$u(2)$ & $0$ & $-u(1)$ & $-u(2)$ & $0$ & $0$ & $H\langle u,v\rangle/2$ & $E\langle u,v\rangle$\\\hline
$v(1)$ & $v(2)$ & $0$ & $v(1)$ & $-F\langle u,v\rangle$ & $-H\langle u,v\rangle/2$ & $0$ & $0$\\\hline 
$v(2)$ & $0$ & $-v(1)$ & $-v(2)$ & $-H\langle u,v\rangle/2$ & $-E\langle u,v\rangle$ & $0$ & $0$\\\hline
\end{tabular}
\end{center}
for $\langle u,v\rangle\in \mathbb{F},$ where $\{u(1), u(2)\}$ and $\{v(1), v(2)\}$ are basis of $V_{2}(u)$ and $V_{2}(v)$, respectively. 
\end{ex}

\subsection{The general case}

We can easily drop the assumption that $mL\neq 0$ for any $0\neq m$ from $\mathcal{M}$.

Let $\mathcal{M}$ be a Malcev algebra in $\mathcal{H}$ containg $L=\mathfrak{sl}_{2}(\mathbb{F}).$ Then, we can consider $\mathcal{M}$ as a Malcev $\mathcal{H}$-module over $L$ and by Corollary $\ref{corol2}$$(i)$, we have 
\begin{center}
    $\mathcal{M}=\textup{Ann}_{\mathcal{M}}L \oplus \widehat{N}_{\mathcal{M}}\oplus J_{\mathcal{M}} 
$,
\end{center}
where $\textup{Ann}_{\mathcal{M}}L=\{m \in \mathcal{M}: m L=0\}$,    
\begin{center}    $\widehat{N}_{\mathcal{M}}=\displaystyle\sum_{0\neq\alpha\in\alpha(\mathcal{M},L,L)} \oplus L\alpha=\displaystyle\sum_{i}\oplus\overline{\mathfrak{sl}_{2}(\mathbb{F})}_{i}$,
\end{center}
and $J_{\mathcal{M}}=\displaystyle\sum_{i} \oplus V_{2i},$ where $V_{2i}$ is a 2-dimensional non-Lie Malcev module for $L.$ 

\begin{lem}\label{lema 3}
$\textup{Ann}_{\mathcal{M}}L$ is a subalgebra of $\mathcal{M}$, $(\textup{Ann}_{\mathcal{M}}L)\widehat{N}_{\mathcal{M}}=0$ and
$(\textup{Ann}_{\mathcal{M}}L)J_{\mathcal{M}}\subseteq J_{\mathcal{M}}.$
\end{lem}
\begin{dem}
Let $m,m'\in \textup{Ann}_{\mathcal{M}}L,$ then 
\begin{equation*}
    (mm')L=(mm')L^{2}=(mm')(LL)\overset{(\ref{ee1})}{=}0.
\end{equation*}
Thus $(\textup{Ann}_{\mathcal{M}}L)^{2}\subseteq \textup{Ann}_{\mathcal{M}}L$ and $\textup{Ann}_{\mathcal{M}}L$ is a subalgebra of $\mathcal{M}$. Furthermore, since $\widehat{N}_{\mathcal{M}}=\displaystyle\sum_{i}\oplus\overline{\mathfrak{sl}_{2}(\mathbb{F})}_{i},$ then 
\begin{equation*}
    (\textup{Ann}_{\mathcal{M}}L)\widehat{N}_{\mathcal{M}}=
(\textup{Ann}_{\mathcal{M}}L)\displaystyle\sum_{i}\oplus\overline{\mathfrak{sl}_{2}(\mathbb{F})}_{i}=0
\end{equation*}
and $(\textup{Ann}_{\mathcal{M}}L)\widehat{N}_{\mathcal{M}}=0$.

Now as $J_{\mathcal{M}}=V(1)\oplus V(2)$, we will prove that $(\textup{Ann}_{\mathcal{M}}L)V_{2}(u)\subseteq J_{\mathcal{M}}$ for any $u\in V=V(1).$
Consider $m\in \textup{Ann}_{\mathcal{M}}L$, then 
\begin{equation*}
\begin{split}
-u(1)mF=u(1)m\cdot FH&\overset{(\ref{ee1})}{=}u(1)Fm H+FmHu(1)+mHu(1) F+Hu(1)Fm\\                 &\overset{(\ref{ee35})}{=}-u(1)Fm=0,
\end{split}
\end{equation*}
we get 
\begin{equation}\label{q1}
    u(1)m F=0.
\end{equation}
Also,
\begin{equation*}
\begin{split}
-u(2)m F=u(2)m \cdot FH&\overset{(\ref{ee1})}{=}u(2)FmH+ FmHu(2)+mHu(2)F+Hu(2)Fm\\
                    &=-u(1)m H+u(2)Fm=-u(1)mH-u(1)m
\end{split}
\end{equation*}
and 
\begin{equation}\label{q2}
    u(2)mF=u(1)mH+u(1)m. 
\end{equation}
Similarly,
\begin{equation*}
\begin{split}
\frac{1}{2}u(1)mH=u(1)m \cdot EF&\overset{(\ref{ee1})}{=}u(1)EmF+EmFu(1)+mFu(1)E+Fu(1)Em=u(2)mF\\
                    &\overset{(\ref{q2})}{=}u(1)mH+u(1)m,
\end{split}
\end{equation*}
then
\begin{equation}\label{q3}
    u(1)mH=-2u(1)m
\end{equation}
and replacing $(\ref{q3})$ in $(\ref{q2}),$ we obtain
\begin{equation}\label{q4}
    u(2)mF=-u(1)m.
\end{equation}
Further,
\begin{equation*}
    u(2)mE=u(2)m\cdot EH\overset{(\ref{ee1})}{=}u(2)EmH+EmHu(2)+mHu(2)E+Hu(2)Em
    =u(2)Em=0;
\end{equation*}
thus 
\begin{equation}\label{q5}
    u(2)mE=0.
\end{equation}
In addition,
\begin{equation*}
\begin{split}
u(1)mE=u(1)m\cdot EH&\overset{(\ref{ee1})}{=}u(1)EmH+EmHu(1)+mHu(1)E+Hu(1)Em\\
                    &=u(2)mH-u(1)Em=u(2)mH-u(2)m,
\end{split}
\end{equation*}
so
\begin{equation}\label{q6}
    u(1)mE=u(2)mH-u(2)m.
\end{equation}
Finally,
\begin{equation*}
\begin{split}
\frac{1}{2}u(2)mH=u(2)m\cdot EF&\overset{(\ref{ee1})}{=}u(2)EmF+EmFu(2)+mFu(2)E+Fu(2)Em\\
&                    =u(1)Em=u(2)m,
\end{split}
\end{equation*}
then 
\begin{equation}\label{q7}
    u(2)mH=2u(2)m
\end{equation}
and replacing $(\ref{q7})$ in $(\ref{q6}),$ we obtain
\begin{equation}\label{q8}
    u(1)mE=u(2)m.
\end{equation}

The equations $(\ref{q1})$, $(\ref{q3})$, $(\ref{q4})$, $(\ref{q5})$, $(\ref{q7})$ and $(\ref{q8})$ give
\begin{center}
    $u(1)mH = -2u(1)m$,\,\,\, $u(1)mF = 0$,\,\,\, $u(1)mE = u(2)m$
\end{center}
\begin{center}
    $u(2)mH=2u(2)m$,\,\,\, $u(2)mF=-u(1)m$,\,\,\, $u(2)mE=0.$
\end{center}
Denote 
$u'(1)=u(1)m$, $u'(2)=u(2)m,$
then 
\begin{center}
    $u'(1)H = -2u'(1)$,\,\,\, $u'(1)F = 0$,\,\,\, $u'(1)E = u'(2),$
\end{center}
\begin{center}
    $u'(2)H=2u'(2)$,\,\,\, $u'(2)F=-u'(1)$,\,\,\, $u'(2)E=0.$
\end{center}
Thus, the set $\{u'(1),u'(2)\}$ form a basis of an $L$-module of type $V_{2}$; so $(\textup{Ann}_{\mathcal{M}}L)V_{2}(u) \subseteq J_{\mathcal{M}}$, then $(\textup{Ann}_{\mathcal{M}}L)J_{\mathcal{M}}\subseteq J_{\mathcal{M}}.$

The lemma is proved.
\end{dem}

\begin{coro}\label{coro3}
The decomposition 
$\mathcal{M}=(\textup{Ann}_{\mathcal{M}}L \oplus \widehat{N}_{\mathcal{M}})\oplus J_{\mathcal{M}}$ is a $\mathbf{Z}_{2}$-graded algebra with 
$\mathcal{M}_{0}=\textup{Ann}_{\mathcal{M}}L \oplus \widehat{N}_{\mathcal{M}}$, $\mathcal{M}_{1}=J_{\mathcal{M}}$. 
\end{coro}
\begin{dem}
It is clear that $\textup{Ann}_{\mathcal{M}}L$ is a Lie $L$-module. Moreover, from Lemma $\ref{l3}$ we know that $J^{2}_{\mathcal{M}}$ is a Lie $L$-module, then $J^{2}_{\mathcal{M}}\subseteq \textup{Ann}_{\mathcal{M}}L \oplus \widehat{N}_{\mathcal{M}}.$ Thus, the result is a consequence of Lemmas $\ref{l2},$ $\ref{lema 3}$ and Corollary \ref{coro1}.

The corollary is proved.
\end{dem}

Furthermore, as already was mentioned, from Lemma $\ref{l3}$ we have $J^{2}_{\mathcal{M}}\subseteq \textup{Ann}_{\mathcal{M}}L \oplus \widehat{N}_{\mathcal{M}}.$ 
A similar process to the proof of Proposition \ref{p3}, leads us to
\begin{equation*}
\begin{split}
u(1)v(1)&=m+F\langle u,v\rangle,\\
u(1)v(2)=u(2)v(1)&=m+\frac{1}{2}H\langle u,v\rangle,\\
u(2)v(2)&=m+E\langle u,v\rangle.\\
\end{split}
\end{equation*}
for some $m\in \textup{Ann}_{\mathcal{M}}L$ and $\langle u,v\rangle\in U,$ where $U$ is a unital commutative associative algebra. But from $(\ref{e6})$, $u(1)v(2)F=\dfrac{1}{2}u(1)v(1),$ then

\begin{center}
    $(m+\dfrac{1}{2}H\langle u,v\rangle)F=\dfrac{1}{2}(m+F\langle u,v\rangle),$
\end{center}
implies $m=0$. Thus, $J^{2}_{\mathcal{M}}\subseteq  \widehat{N}_{\mathcal{M}};$ so by Lemma \ref{l2} and Corollary \ref{coro1}, $\widehat{N}_{\mathcal{M}}\oplus J_{\mathcal{M}}=\mathfrak{sl}_{2}(U) \oplus V(1) \oplus V(2)$ is a subalgebra of $\mathcal{M}.$ 

As $\widehat{N}_{\mathcal{M}}\oplus J_{\mathcal{M}}$ is a subalgebra, we can apply to $\widehat{N}_{\mathcal{M}}\oplus J_{\mathcal{M}}$ the results obtained in subsection \ref{subsection}. Thus, we can write 
\begin{center}
    $\mathcal{M}=\textup{Ann}_{\mathcal{M}}L \oplus (\mathfrak{sl}_{2}(U) \oplus V^2),$
\end{center}
where $\widehat{N}_{\mathcal{M}}\cong \mathfrak{sl}_{2}(U)$, $J_{\mathcal{M}}=V^{2}$ and 
$U$ satisfies
\begin{center}
    $U\subseteq \Gamma(\mathfrak{sl}_{2}(U) \oplus V^2)$
\end{center}
and $V$ is a commutative $U$-bimodule. Also there exists a $U$-bilinear skew-symmetric mapping $\langle\cdot , \cdot\rangle: V \times V \rightarrow U$ such that $\langle V, V\rangle \subseteq U$ and formula $(\ref{ee28})$ holds for any $u, v, w \in V$. 
Moreover in $\mathfrak{sl}_{2}(U) \oplus V^2$ is defined the product $(\ref{ee33}).$

Therefore, we have our main results of this section.

\begin{pro}
The Malcev algebra $\mathcal{M}=\textup{Ann}_{\mathcal{M}}L \oplus (\mathfrak{sl}_{2}(U) \oplus V^2)$ in $\mathcal{H}$ containing $\mathfrak{sl}_{2}(\mathbb{F})$ is isomorphic to  $\textup{Ann}_{\mathcal{M}}L \oplus \mathcal{M}_{7}(U)$ if and only if  there exist $\langle u, v\rangle=1$ for $u=(\frac{1}{2},0), v=(0,-\frac{1}{2}) \in V.$
\end{pro}
\begin{dem}
From Proposition \ref{p2}, the subalgebra $\mathfrak{sl}_{2}(U) \oplus V^2$  is isomorphic
to the exceptional Malcev algebra $\mathcal{M}_{7}(U)$ if and only if  there exist $\langle u, v\rangle=1$ for $u=(\frac{1}{2},0), v=(0,-\frac{1}{2}) \in V.$ 

The proposition is proved. 
\end{dem}

Finally from the Cayley–Dickson process, we get the following result.

\begin{theo}
The Malcev algebra $\mathcal{M}=\textup{Ann}_{\mathcal{M}}L \oplus (\mathfrak{sl}_{2}(U) \oplus V^2)$ in $\mathcal{H}$ containing $\mathfrak{sl}_{2}(\mathbb{F})$ is isomorphic to  $\textup{Ann}_{\mathcal{M}}L \oplus\widetilde{\mathcal{M}}\left(\mathfrak{sl}_{2}(U), \alpha\right)$, where $\alpha=-\langle(\frac{1}{2},0),(0,\frac{1}{2})\rangle$. 
\end{theo}
\begin{dem}
By Theorem \ref{t2}, the subalgebra $\mathfrak{sl}_{2}(U) \oplus V^2$ is isomorphic to the Malcev algebra $\widetilde{\mathcal{M}}\left(\mathfrak{sl}_{2}(U), \alpha\right),$ where $\alpha=-\langle(\frac{1}{2},0),(0,\frac{1}{2})\rangle$. 

The theorem is proved.
\end{dem}

\section{Funding}

The article was supported by Resolution of the Organizing Commission $\textup{N}^{\circ}$ $629-2022$-UNAB of the National University of Barranca. 

Also the  author  gratefully  acknowledges  financial  support by CONCYTEC-PROCIENCIA within the framework of the call “Proyecto Investigación Básica 2019-01” 
[380-2019-FONDECYT].

\section{Statements and Declarations}

\subsection{Competing Interests} The author states that there is no Conflict of interest.

%



\end{document}